\documentstyle{amsppt}
\magnification 1000 \hsize=6.1in \vsize=8.75in
\input amstex
\TagsOnRight

\document
 
\def\s#1#2{\langle \,#1 , #2 \,\rangle}

\def\H{{\frak H}}
\def\F{{\frak F}}
\def\C{{\Bbb C}}
\def\R{{\Bbb R}}
\def\Z{{\Bbb Z}}

\def\G{\Gamma}
\def\g{\gamma}
\def\L{\Lambda}
\def\ee{\varepsilon}

\def\ca{{\frak a}}
\def\cb{{\frak b}}

\def\ci{{\infty}}

\def\sa{{\sigma_\frak a}}
\def\sb{{\sigma_\frak b}}

\def\p{\endproclaim \flushpar {\bf Proof: }}

\def\bs{\ \ $\qed$ \vskip 3mm}

\def\Im{{\text{\rm Im}}}
\def\Re{{\text{\rm Re}}}

\topmatter
\title{The dimensions of spaces of holomorphic second-order automorphic 
forms and their cohomology}
\endtitle
\author{NIKOLAOS DIAMANTIS (University of Nottingham) \\
CORMAC O'SULLIVAN (City University of New York)}
\endauthor
\endtopmatter

\NoRunningHeads

$$\text{\bf 1. Introduction}$$

In the paper \cite{KZ}, Kleban and Zagier find that the 
study of crossing probabilities in percolation theory leads naturally to 
holomorphic second-order modular forms.
In this paper we answer a question posed by Zagier and compute 
the exact dimensions of these second-order spaces for even weight. 
We also establish a cohomological interpretation
of these spaces which is
analogous to that of Eichler and Shimura for usual modular forms.

\vskip .10in Second-order automorphic forms have recently arisen in several other 
contexts besides
percolation theory, for example Eisenstein series with modular symbols in 
the papers \cite{G},
\cite{O1}, \cite{PR} and GL$(2)$ converse theorems in \cite{F} and \cite{FW}.
In \cite{CDO} classification theorems for smooth second-order forms were provided, but we were only able to prove 
an upper bound for the dimensions of the holomorphic spaces.
Here we establish a formula
for the dimension confirming that it equals the upper bound given in
\cite{CDO} when the weight is not 2. The weight 2 case was the hardest to resolve, requiring the analytic continuation of a number of related series. Interestingly, in this case the dimension differs from the natural upper bound by 1. These results may have 
consequences for some percolation
theory problems and the converse theorems investigated by Farmer and Wilson.

\vskip .10in We use the knowledge of the dimensions to express the second-order spaces 
in terms of Eichler cohomology. As
in the classical theorem of Eichler and Shimura, the isomorphism is given
explicitly by period polynomials thereby maintaining the connection
of such maps with values of (twisted) $L$-functions.
We expect this isomorphism will yield applications analogous to those of
the Eichler-Shimura isomorphism. 
We also hope it
will help us in obtaining a natural geometric interpretation for 
second-order automorphic forms.

$$\text{\bf 2. Definitions and Statement of Main Results}$$

Let $\G \subset PSL_2(\R)$ be a 
Fuchsian group of the first kind acting on the upper half plane $\H$ 
with non compact quotient $\G\backslash\H$.  As usual we write $x+iy=z\in 
\H$. Let $d\mu z$ be 
the hyperbolic volume form $dx dy/y^2$ and $V$ the volume of 
$\G\backslash\H$. For a fundamental domain $\F$ fix representatives of 
the inequivalent cusps in $\overline \F$ and give them labels such as 
$\ca,\cb$. Use the corresponding scaling matrices $\sa,\sb$ to give 
convenient local coordinates near these cusps as in \cite{I1}, Chapter 2 
for example. The subgroup $\G_\ca$ is the set of elements of $\G$ fixing 
$\ca$ and 
$$
\sa^{-1} \G_\ca \sa= \G_\ci=
\left\{ \pm \left(\smallmatrix 1 & m \\ 0 & 1
\endsmallmatrix\right)
\; \big | \; \ m\in {\Bbb Z}\right\}.
$$
The slash operator $|_k$ defines an action of $PSL_2(\R)$ on 
functions $f:\H \mapsto \C$ by
$$
(f|_k \gamma)(z)=f(\gamma z) (c z+d)^{-k} 
$$
with $\gamma= \left ( \smallmatrix  * & * \\  c & d
 \endsmallmatrix \right )$ in $PSL_2(\R)$.
Extend the action to $\Bbb C [PSL_2(\Bbb R)]$ by linearity. In this context we set $j(\g ,z)= cz+d$ for later use.

\vskip .10in We shall also require the generators of the group $\G$. Suppose $\Gamma \backslash \H$ has genus $g,$ $r$ elliptic fixed 
points
and $p$ cusps, then there are $2g$ hyperbolic elements $\gamma_i,$
$r$ elliptic elements $\epsilon_i$ and $p$ parabolic elements 
$\pi_i$ generating $\Gamma$ and satisfying the $r+1$ 
relations: 
$$[\gamma_1, \gamma_{g+1}]\dots
[\gamma_{g}, \gamma_{2g}]\epsilon_1 \dots \epsilon_r \pi_1 \dots 
\pi_{p}=1, \qquad \epsilon_j^{e_j}=1 \tag 2.1$$
for $1 \leqslant j \leqslant r$ and integers $e_j \geqslant 2.$ 
Here $[a, b]$ denotes the
commutator $aba^{-1}b^{-1}$ of $a, b.$ (cf. [I1] Proposition 2.6).

\vskip .10in Let $k$ be an integer. 
As explained below, a first or second-order form $f$ is defined from the following list of 
alternatives. It may satisfy either
\roster
\item"{\bf H.}" $f:\H \rightarrow \C$ is holomorphic or
\item"{\bf S.}" $f:\H \rightarrow \C$ is smooth.
\endroster
Its automorphy condition may involve, for some vector space of functions $V$,
\roster
\item"{\bf A(}V{\bf ).}" $f|_k(\g -1)\in V$ for all  $\g$ in $\G$ or
\item"{\bf P.}" $f|_k(\pi -1)=0$ for all parabolic $\pi$ in $\G$.
\endroster
 Finally, we need to impose growth conditions which may be cuspidal ({\bf C}), or non-cuspidal ({\bf N}):
\roster
\item"{\bf C.}" for each cusp $\ca$, $(f|_k \sa)(z) \ll e^{-c y}$ as $y\rightarrow \ci$ uniformly in $x$ for some constant $c>0$ or
\item"{\bf N.}" for each cusp $\ca$, $(f|_k \sa)(z) \ll y^c$ as $y\rightarrow \ci$ uniformly in $x$ for some constant $c$.
\endroster

\vskip .10in Define $S_k(\G)$, the space of holomorphic, weight $k$, cusp forms for $\G$, to be the $\C$-vector space of  functions $f$ such that {\bf H}, {\bf A(}0{\bf )} and {\bf C} hold. 
The corresponding space of modular forms, $M_k(\G)$, satisfies {\bf H}, {\bf A(}0{\bf )} and {\bf N}.
Both these spaces are always finite dimensional. Using the Riemann-Roch theorem, exact formulas for their dimensions are given in \cite{Sh}, Theorems 2.23, 2.24 (see Section 4). Following \cite{DKMO}, the space $S_k^2(\G)$ of holomorphic, weight $k$, second-order cusp forms consists of  functions satisfying {\bf H}, {\bf A(}$S_k(\G)${\bf )}, {\bf P} and {\bf C}.
This new space is similar to $S_k(\G)$, the only difference being the weaker automorphy rule.
Naturally, we define the space $M_k^2(\G)$ with the conditions {\bf H}, {\bf A(}$M_k(\G)${\bf )}, {\bf P} and {\bf N}.
Lastly, replacing {\bf H} with {\bf S}, we obtain the smooth spaces $\tilde S_k(\G)$,  $\tilde M_k(\G)$, $\tilde S_k^2(\G)$ and $\tilde M_k^2(\G)$ (with {\bf A(}$\tilde S_k(\G)${\bf )} and {\bf A(}$\tilde M_k(\G)${\bf )} in the definitions of these last two). They are infinite dimensional in general. Restricting to certain eigenspaces of the Laplacian does yield finite dimensional spaces, see \cite{CDO}.

\vskip .10in To summarize, the main holomorphic spaces used subsequently are:
$$
\matrix
\bullet &  S_k(\G) & \text{defined with} & \text{{\bf H}}, & \text{{\bf A(}}0\text{{\bf )}} & \text{and {\bf C},} & \qquad \qquad \qquad \qquad \qquad \qquad \qquad \\
\bullet &  S_k^2(\G) & \text{defined with} & \text{{\bf H}}, & \text{{\bf A(}}S_k(\G)\text{{\bf ), P}} & \text{and {\bf C},} & \qquad \qquad \qquad \qquad \qquad \qquad \qquad \\
\bullet &  R_k^2(\G) & \text{defined with} & \text{{\bf H}}, & \text{{\bf A(}}S_k(\G)\text{{\bf ), P}} & \text{and {\bf N},} & \qquad \qquad \qquad \qquad \qquad \qquad \qquad \\
\bullet &  M_k(\G) & \text{defined with} & \text{{\bf H}}, & \text{{\bf A(}}0\text{{\bf )}} & \text{and {\bf N},} & \qquad \qquad \qquad \qquad \qquad \qquad \qquad \\
\bullet &  M_k^2(\G) & \text{defined with} & \text{{\bf H}}, & \text{{\bf A(}}M_k(\G)\text{{\bf ), P}} & \text{and {\bf N}.} & \qquad \qquad \qquad \qquad \qquad \qquad \qquad  
\endmatrix
$$
We included the space $R_k^2(\G)$ above which has some interesting elements as we shall see in Section 5. Their smooth counterparts are:
$$
\matrix
\bullet &  \tilde S_k(\G) & \text{defined with} & \text{{\bf S}}, & \text{{\bf A(}}0\text{{\bf )}} & \text{and {\bf C},} & \qquad \qquad \qquad \qquad \qquad \qquad \qquad \\
\bullet &  \tilde S_k^2(\G) & \text{defined with} & \text{{\bf S}}, & \text{{\bf A(}}\tilde S_k(\G)\text{{\bf ), P}} & \text{and {\bf C},} & \qquad \qquad \qquad \qquad \qquad \qquad \qquad \\
\bullet &  \tilde M_k(\G) & \text{defined with} & \text{{\bf S}}, & \text{{\bf A(}}0\text{{\bf )}} & \text{and {\bf N},} & \qquad \qquad \qquad \qquad \qquad \qquad \qquad \\
\bullet &  \tilde M_k^2(\G) & \text{defined with} & \text{{\bf S}}, & \text{{\bf A(}}\tilde M_k(\G)\text{{\bf ), P}} & \text{and {\bf N}.} & \qquad \qquad \qquad \qquad \qquad \qquad \qquad 
\endmatrix
$$
In the paper \cite{DKMO} we examine the effects of weakening or altering the conditions defining each of these spaces. Here we prove three main results.  

\proclaim{Theorem 2.1} For $k$ in $2\Z$ and $\G\backslash\H$ non compact with genus $g$ we have
$$
\align
\dim S^2_{k}(\G) & = 0 \text{ if } k \leqslant 0,\\ 
\dim S^2_2(\G) & =\cases 0 \text{ if }\dim S_2(\G)=0, \\ (2g+1)\dim S_2(\G) -1\text{ otherwise,} \endcases\\ 
\dim S^2_k(\G) & =(2g+1)\dim S_k(\G) \text{ if }k \geqslant 4.
\endalign
$$
\endproclaim
\proclaim{Theorem 2.2} For $k$ in $2\Z$ and $\G\backslash\H$ non compact with genus $g$ we have
$$
\align
\dim M^2_{k}(\G)  &= 0 \text{ if } k \leqslant -2,\\
\dim M^2_{0}(\G)  &= g+1,\\  
\dim M^2_k(\G)  &=(2g+1)\dim M_k(\G) \text{ if }k \geqslant 2.
\\
\endalign
$$
\endproclaim

\flushpar

\proclaim{Theorem 2.3}  We have the 
isomorphisms
$$
\align
\frac{M_2^2(\Gamma)}{M_2(\Gamma)} & \oplus \frac{\bar S_2^2(\Gamma)}{\bar
S_2(\Gamma)} \oplus \Bbb C  \cong H^1_{!}(\Gamma, \Bbb C ),\\
\frac{M_k^2(\Gamma)}{M_k(\Gamma)} & \oplus \frac{\bar S_k^2(\Gamma)}{\bar
S_k(\Gamma)}  \cong H^1_{!}(\Gamma, P_{k-2} ) \text{ \ for even }k 
\geqslant 4
\endalign
$$
where $P_{k-2}$ is the space of all polynomials of degree at most $k-2$ 
and $H^1_{!}(\Gamma, \C )$, $H^1_{!}(\Gamma, P_{k-2} )$  are associated 
cohomology groups as defined in Section 7.
\endproclaim

We note that Goldfeld considered certain subspaces of $S^2_k(\G)$ and $M^2_k(\G)$ in \cite{G} (preceeding the work in \cite{KZ}) and calculated their dimensions. See also \cite{DO}, Proposition 15, for a further generalization.

\vskip .10in The proofs of Theorems 2.1 and 2.2 are constructive in that, in each case, a basis of second-order Poincar\'e series is given. In the weight 2 case this relies on the analytic continuation of various series. An example of the type of results we prove is the following.
For $m \in \Z$ and non-negative define the non-holomorphic, weight $k$
Poincar\'e series
$$
U_{\ca m}(z,s,k)=\sum_{\g \in \G_\ca \backslash \G}   \Im(\sa^{-1} \g 
z)^s e(m \sa^{-1} \g z) \ee(\sa^{-1} \g, z)^{-k} \tag 2.2
$$
for $\Re (s)>1$ where $\ee(\g, z)$ is the `non-holomorphic' weight factor
$j(\g,z)/|j(\g,z)|$ and, as usual, $e(z)=e^{2\pi i z}$.

\proclaim{Proposition 2.4} For $k \in 2\Z$, $m > 0$ and some $\delta>0$ the Poincar\'e series
$U_{\ca m}(z,s,k)$  has a
continuation to an analytic function for all $s$ with $\Re (s) >1-\delta$. We
have
$$
U_{\ca m}(z,s,k) \ll   y_\G(z)^{1/2}
$$
for these $s$ with the implied constant depending on $s$, $m$, $k$ and $\G$.
\endproclaim
\flushpar
Here $y_\G(z)$ is the invariant height
$$
y_\G(z)=\max_\ca (\max_{\g \in \G}(\Im( \sa^{-1}\g z))).
$$

\newpage

$$\text{\bf 3. Upper bounds for $\dim S^2_k(\G)$ and $\dim M^2_k(\G)$}$$

Let Hom$_0(\G, \C)$ be the space of homomorphisms from $\G$ to $\C$ that are zero  on parabolic elements of $\G$. According to a special case of the Eichler-Shimura isomorphism (see Section 7 for a complete statement),
for any such homomorphism $L$ there exist unique $l_1$, $l_2$ in $S_2(\G)$ so that
$$
L(\g)=\int_z^{\g z} l_1(w)\,dw+\overline{\int_z^{\g z} l_2(w)\,dw} \tag 3.1
$$
for all $\g$ in $\G$. We define the {\it modular symbol}
$$
\s{\g}{l}=\int_z^{\g z} l(w)\,dw \tag 3.2
$$
for $\g$ in $\G$ and $l$ in $S_2(\G).$
The value of $\s{\g}{l}$ is independent of the choice of base 
point $z$ in $\H$. 
For $1\leqslant i \leqslant 2g$ we may define
$$
\L_i(z)=\int_{z_0}^z \lambda_i(w)\, dw + \overline{\int_{z_0}^z \mu_i(w)\, dw} \tag 3.3
$$
for certain $\lambda_i$, $\mu_i$ in $S_2(\G)$ and $z_0$ an arbitrary 
fixed element of $\H$ (usually taken to be the imaginary number $i$ or $i\ci$) to satisfy
$$
\L_i(\g_j z)-\L_i( z)=\delta_{ij}
$$
for the hyperbolic generators and also 
$$
\L_i(\g z)-\L_i( z)=0
$$
for the other parabolic and elliptic generators $\g$ of $\G$. Thus $\{L_i\}_{i=1}^{2g}$ defined by
$$L_i(\g)= \L_i(\g z)-\L_i( z) \tag 3.4$$
forms a natural basis for Hom$_0(\G, \C)$ dual to our choice of generators for $\G$.

\vskip .10in From the map
$$
f\mapsto (f|_k(\g_1 -1), \dots ,f|_k(\g_{2g} -1))=(f_1, \dots ,f_{2g}),
$$
it may be shown, as in \cite{DKMO} Theorem 4, that any element $f$ of $S_k^2(\G)$ must satisfy
$$
f(z)=\sum_{i=1}^{2g}f_i(z)\L_i(z)+\phi(z)\tag 3.5
$$
where each $f_i=f|_k(\g_{i} -1)$ is in $S_k(\G)$ and $\phi$ is in $\tilde{S}_k(\G)$. 
Apply $\frac{d}{d\overline{z}}$ to both sides of $(3.5)$ to get
$$
\frac{d}{d\overline{z}}\phi=-\sum_{i=1}^{2g}f_i \overline{\mu_i}.\tag 3.6
$$
So $\phi$, if it exists, is uniquely specified by $(3.6)$ up to addition of an element of $S_k(\G)$. Counting the degrees of freedom in $(3.5)$ shows that, for $k \in \Z$, 
$$
\dim S_k^2(\G) \leqslant (2g+1)\dim S_k(\G).\tag 3.7
$$
Similarly, $f$ in $M_k^2(\G)$ satisfies $(3.5)$ with $f_i$ in $M_k(\G)$ and $\phi$ in $\tilde{M}_k(\G)$. Thus
$$
\dim M_k^2(\G) \leqslant (2g+1)\dim M_k(\G).\tag 3.8
$$

\vskip .10in In fact the dimensions of these spaces, for weight $k \neq 2$, are given by the upper bounds in $(3.7)$ and $(3.8)$. We will show this in the next section using an extension of the usual Poincar\'e series construction to build linearly independent second-order forms. The case of weight $k=2$ must be singled out for special attention since the Poincar\'e series we need are no longer absolutely convergent. An analysis in Sections 5 and 6 shows that our constructions come up one short of the upper bounds at this weight.

\vskip .10in Before discussing Poincar\'e series, we note that 
there is a much easier way to find elements of $S_k^2(\G)$ and $M_k^2(\G)$. Certainly we have the subspaces $S_k(\G) \subset S_k^2(\G)$ and 
$M_k(\G) \subset M_k^2(\G)$. Also, it may be checked (see Lemma 4.1 in the next section) that any product
$$
f(z)\int_{z_0}^z h(w)\, dw \tag 3.9
$$
for $f$ in $S_k(\G)$ and $h$ in $S_2(\G)$ gives an element of $S_k^2(\G)$. If $f$ is in $M_k(\G)$ then $(3.9)$ yields an element of $M_k^2(\G)$. 
This suggests the natural decomposition:
$$
\align
S_k^2(\G)&=S_k^{2}(\G)^+ \oplus S_k(\G) \oplus S_k^2(\G)^-, \tag 3.10 \\
M_k^2(\G)&=M_k^{2}(\G)^+ \oplus M_k(\G) \oplus M_k^2(\G)^-, \tag 3.11
\endalign
$$
where the spaces $S_k^2(\G)^+$, 
$M_k^{2}(\G)^+$ consist of all linear combinations of elements of the  
form $(3.9)$.
In this way, for $k \geqslant 0$, we get the simple lower bounds
$$
\align
(g+1)\dim S_k(\G)=\dim (S_k^{2}(\G)^+ \oplus S_k(\G)) &\leqslant \dim 
S_k^2(\G),\\
(g+1)\dim M_k(\G)=\dim (M_k^{2}(\G)^+ \oplus M_k(\G)) &\leqslant \dim 
M_k^2(\G).
\endalign
$$

\vskip .10in A second-order form similar to $(3.9)$ appears in the work 
of Kleban and Zagier, \cite{KZ}, on percolation theory. For the Dedekind 
eta function $\eta(z)$ we have that $\eta(z)^4 \in S_2(\G_\theta, 
\chi)$ and that
$\eta(z/2)^8\eta(2z)^8\eta(z)^{-12}$ is an automorphic form with character  
$\overline{\chi}$ for $\G_\theta$ which is holomorphic at $\frak H$, 
vanishing at $\infty$ and has a pole at $1$.
Here $\G_\theta$ is the subgroup of $PSL_2(\Z)$ generated by $\pm\left ( 
\smallmatrix  1 & 2 \\  0 & 1
 \endsmallmatrix \right )$, $\pm\left ( \smallmatrix  0 & -1 \\  1 & 0
 \endsmallmatrix \right )$. Also $\chi$ is a certain character for 
$\G_\theta$. They show in equations (16), (19) of \cite{KZ} that the 
derivative of a certain crossing probability for a rectangle with aspect 
ratio $r$ is 
$$
K(z)=\frac{-16\pi i}{\sqrt{3}}\eta(z)^4 \int_{i\infty}^z \eta(w/2)^8\eta(2w)^8\eta(w)^{-12}\, dw  \tag 3.12
$$
at $z=ir$. Consult \cite{KZ} for details. The expression $K(z)$ is 
a type of second-order form satisfying $K|_2(\g -1) \in S_2(\G_\theta, \chi)$ for all $\g$ in $\G_\theta$, and conditions {\bf H} and {\bf C}.

$$\text{\bf 4. Exact dimensions of $S^2_k(\G)$ and $M^2_k(\G)$ for $k>2$}$$

First we look at the spaces $S_k(\G)$ and $M_k(\G)$ with $k \in 2\Z$.  
Recall that $g$ is the genus  of $\Gamma \backslash \H$ and $p$ the 
number of inequivalent cusps. By Theorems 2.23 and 2.24 of \cite{Sh} we 
have the following:
$$
\align
\dim S_{k}(\G) &=  \dim M_{k}(\G)=0 \text{ for } k<0,\\
\dim S_{0}(\G) &= \cases 1 \text{ if } p=0 \\ 0 \text{ if } p>0 \endcases, \\
\dim S_{2}(\G) &= g,\\
\dim M_{0}(\G) &= 1,\\
\dim M_{2}(\G) &= \cases g \text{ if } p=0 \\ g+p-1 \text{ if } p>0 \endcases,\\
\dim M_{k}(\G) &= \dim S_{k}(\G)+p \text{ if $k\geqslant 4$}.
\endalign
$$
Also for $k\geqslant 4$, with $(2.1)$,
$$
\dim S_k(\G)=(k-1)(g-1)+(k/2-1)p+\sum_{j=1}^r \frac{k(e_j -1)}{2 e_j}.
$$
The analogous results for $k$ odd, $ \neq 1$ also appear in \cite{Sh} Theorem 2.25.

\vskip .10in What are the elements of these spaces?  For $k=0$ the elements of 
$M_0(\G)$ (and $S_0(\G)$ if there are no cusps) are just the constant 
functions. When $k \geqslant 4$ is even and $p>0$, $S_k(\G)$ is spanned by the 
Poincar\'e series 
$$
P_{\ca m}(z)=P_{\ca m}(z)_k=\sum_{\g \in \G_\ca \backslash \G} j(\sa^{-1} 
\g, z)^{-k} e(m \sa^{-1} \g z)\tag 4.1
$$
with $m>0$. 
For example we may fix a single cusp $\ca$ and find a linearly 
independent basis with $\dim S_k(\G)$ different values of $m>0$, see 
\cite{I2, Corollary 3.5}. The extra $p$ dimensions of $M_k(\G)$ come from 
the linearly independent $P_{\ca 0}(z)_k$ as $\ca$ varies over the $p$ 
inequivalent cusps as in \cite{Sa}, Section 1.4. When $m=0$ these series 
$(4.1)$ are called the holomorphic Eisenstein series.
If we let $E_k(\G)$ denote the space of these Eisenstein series then we have the direct sum
$$
M_k(\G)=E_k(\G) \oplus S_k(\G). \tag 4.2
$$
Thus, varying $\ca$ and $m \geqslant 0$ in $(4.1)$ and taking linear 
combinations produces all elements of $M_k(\G)$ for $k \geqslant 4$. We 
describe what happens at weight $k=2$ later in Section 5.

\vskip .10in We will now prove that the dimensions of $S_k^2(\G)$ and 
$M_k^2(\G)$ attain the upper bounds (3.7) and (3.8)  by 
producing enough second-order forms. In some cases this can be achieved 
using results of \cite{Gu} and \cite{K}. However, these techniques would 
not yield explicit bases. To construct explicit bases we use an extension 
of the above  Poincar\'e series. This idea, in different guise, appears
in \cite{G}. Set
$$
P_{\ca m}(z,L)=
P_{\ca m}(z,L)_k=\sum_{\g \in \G_\ca \backslash \G} L(\g) j(\sa^{-1} \g, 
z)^{-k} e(m \sa^{-1} \g z), \tag 4.3
$$
for $m\geqslant 0$ and $L$ in Hom$_0(\G, \C)$. 
To show that this is absolutely 
convergent and holomorphic for $k \geqslant 4$ we need the next lemma. 

\proclaim{Lemma 4.1} For any $f$ in $S_2(\G)$, all $z \in \H$ and any 
cusp $\ca$,
$$
\int_{z_0}^z f(w)\, dw \ll |\log(\Im (\sa^{-1} z))| +1
$$
with an implied constant independent of $z$.
\endproclaim

The proof uses that $y|f(z)| \ll 1$ for any weight 2 cusp form and also 
that $f(z+1)=f(z)$. See \cite{DKMO}, Lemma 3 for details.

\proclaim{Proposition 4.2} For $L$ in $\operatorname{Hom}_0(\G, \C)$ and $k \geqslant 4$ it is the case that 
$$
\align
P_{\ca m}(z,L)_k & \in M_k^2(\G) \text{ \ \ if \ } m=0,\\
P_{\ca m}(z,L)_k & \in S_k^2(\G) \text{ \ \ if \ } m>0.
\endalign
$$
\p
For $L$ as  above, all $z \in \H$ and any cusp $\ca$ it follows from 
$(3.1)$ and Lemma 4.1 that
$$
L(\g) \ll |\log(\Im (\sa^{-1}\g z))|+|\log(\Im (\sa^{-1} z))| +1
$$
with an implied constant independent of $z$ and $\g \in \G$.
It is also true that $|\log y| <y+1/y$ and hence $|\log y| 
<\varepsilon^{-1} (y^\varepsilon+y^{-\varepsilon})$ for any $y>0$ and 
$\varepsilon>0$. Therefore 
$$
\align
P_{\ca m}(z,L)_k &\ll \sum_{\g \in \G_\ca \backslash \G} (\Im(\sa^{-1} \g 
z)^{\varepsilon}+\Im(\sa^{-1} \g z)^{-\varepsilon} + \Im(\sa^{-1} 
z)^{\varepsilon}+\Im(\sa^{-1} z)^{-\varepsilon}+1) \\
& \qquad \qquad \times |j(\sa^{-1} \g, z)|^{-k}\\
&= y^{-k/2}\sum_{\g \in \G_\ca \backslash \G}  (\Im(\sa^{-1} \g 
z)^{k/2+\varepsilon}+\Im(\sa^{-1} \g z)^{k/2-\varepsilon}) \\
& \qquad \qquad + y^{-k/2}(\Im(\sa^{-1} 
z)^{\varepsilon}+\Im(\sa^{-1} z)^{-\varepsilon}+1) \sum_{\g \in \G_\ca \backslash \G}  \Im(\sa^{-1} \g 
z)^{k/2} \tag 4.4
\endalign
$$
for any $\varepsilon>0$. (The implied constant depends on 
$\varepsilon.$)
Now the usual non-holomorphic Eisenstein series
$$
E_{\ca}(z,s)=\sum_{\g \in \G_\ca \backslash \G}   \Im(\sa^{-1} \g z)^s 
,\tag 4.5
$$
is known to be absolutely convergent for $s$ with $\Re(s)>1$ (and uniformly convergent for $s$ in compact sets there), see \cite{I1} Chapter 3 and \cite{Sa} Section 1.4. Moreover it has the Fourier expansion at the cusp $\cb$
$$
\align
E_{\ca}(\sb z,s)&=\delta_{\ca \cb}y^s +\phi_{\ca \cb}(s)y^{1-s}+\sum_{m \neq 0} \phi_{\ca \cb}(m,s)W_s(mz) \tag 4.6\\
&=\delta_{\ca \cb}y^s +\phi_{\ca \cb}(s)y^{1-s}+O(e^{-2\pi y}) \tag 4.7
\endalign
$$
as $y \rightarrow \infty$ with an implied constant depending only on $s$ 
and $\G$. This is \cite{I1}, (6.20) and is in fact valid for all $s$ in $\C$. The Whittaker function $W_s(z)$ is described in \cite{I1}, (1.26).

\vskip .10in By comparing (4.4)  with (4.5) we require $k/2-\varepsilon>1$ and hence 
$k>2$ for the absolute and uniform convergence of $P_{\ca m}(z,L)_k$. 
Therefore, for $k \ge 4,$ the series $P_{\ca m}(z,L)_k$ are absolutely and uniformly convergent (for $z$ in compact sets in $\H$, say) and  
satisfy {\bf H}.

\vskip .10in We also easily have 
$$
P_{\ca m}(\g z,L)_k j(\g,z)^{-k}=P_{\ca m}( z,L)_k +L(\g^{-1})P_{\ca m}( z)_k
\tag 4.8$$
which implies that $P_{\ca m}(z,L)_k$ satisfies {\bf A(}$M_k(\G)${\bf )} for $m=0$ and 
{\bf A(}$S_k(\G)${\bf )} for $m>0$.

\vskip .10in Finally, we verify that our functions satisfy the cuspidal growth 
condition {\bf C} by considering the Fourier 
expansion of these series at $\cb$. Since we have established {\bf H} we must have 
$$
j(\sb, z)^{-k}P_{\ca m}(\sb z,L)_k=\sum_{n \in \Z} a_{\ca \cb}(n) e(nz)
$$
for some constants $a_{\ca \cb}(n)$. Also (noting that $L(I)=0$ for $I$ 
the identity element of $\G$) and taking into account (4.7) we have
$$
\align
j(\sb, z)^{-k}P_{\ca m}(\sb z,L)_k &=\sum_{\g \in \G_\ca \backslash \G} 
L(\g) j(\sa^{-1} \g \sb, z)^{-k} e(m \sa^{-1} \g \sb z)\\
&\ll y^{-k/2} \sum_{\g \in \G_\ca \backslash \G, \g \neq \G_\ca} L(\g) \Im(\sa^{-1} \g \sb z)^{k/2}\\
&\ll y^{-k/2} \left| E_{\ca}(\sb z,k/2-\epsilon)-\delta_{\ca \cb}y^{k/2-\epsilon}\right|\\
&\ll 
y^{1-k+\epsilon}.
\endalign
$$
Therefore, for $m \geqslant 0$, we find that $a_{\ca \cb}(n)=0$ for $n 
\leqslant 0$. Consequently 
$P_{\ca m}(z,L)_k$ satisfies the cuspidal growth condition {\bf C}. This 
is perhaps surprising in the case where $m=0$. It means that on the 
fundamental domain $\F$ (corresponding to $\sa$, $\sb$ etc), $P_{\ca 
0}(z,L)_k$ has exponential decay at its cusps. This will not be the case 
on any translates $\g \F$ of $\F$ with $L(\g) \neq 0$ by the automorphy 
condition. \bs

\proclaim{Theorem 4.3}
For $k\geqslant 4$ and $g$ the genus of $\G\backslash\H$ we have
$$
\align
 \dim S_k^2(\G) &= (2g+1)\dim S_k(\G), \tag 4.9\\
\dim M_k^2(\G) &= (2g+1)\dim M_k(\G). \tag 4.10
\endalign
$$
\p
Note that, for a fixed cusp $\ca$, the Poincar\'e series $P_{\ca 
j}(z,L_i)$ are all linearly independent, by (4.8), as $j>0$ runs over 
integers  yielding a basis $P_{\ca j}(z)$ for $S_k(\G)$ and as $i$ 
runs over integers in $\{1, \dots, 2g\}$ yielding a basis $L_i$ of 
$\operatorname{Hom}_0(\G,\C)$. These series are also linearly independent 
of $S_k(\G)$. In this way we obtain $(4.9)$.

\vskip .10in 
A similar argument, using the fact that $P_{\ca 0}(z)$ with $\ca$ running over the inequivalent cusps of $\G \backslash \H$ form a basis for $E_k(\G)$, yields $(4.10)$. \bs

This result also clarifies the direct sums $(3.10)$ and $(3.11)$. For example, a second-order cusp form is in the space $S_k^2(\G)^+$ if and only if it is a linear combination of Poincar\'e series $P_{\ca m}(z,L)$ with $L$ in the subspace of Hom$_0(\G, \C)$ generated by the modular symbols $(3.2)$.
The space $S_k^2(\G)^-$ consists of linear combinations of Poincar\'e series $P_{\ca m}(z,\overline{L})$, with the conjugates of the modular symbols.

$$\text{\bf 5. Calculating $\dim S^2_2(\G)$}$$

We first recall some definitions and terminology and put in 
place the general framework used in \cite{CO}, \cite{JO}. We will then 
specialize to the objects we need.

\vskip .10in For $f_1$, $f_2$ in $S_k(\G)$ the Petersson inner product is defined by
$$
\s{f_1}{f_2}=\int_{\G \backslash \H} f_1(z) \overline{f_2(z)} y^k \, d\mu z.
$$
The weight, $k$, of the inner product should be clear from the context; in the 
rest of this section they are weight 2, in Section 8 they are all weight 
0. (Also recall that we are using this notation for the modular 
symbols (3.2).)

\vskip .10in Following \cite{I1}, (2.42) we recall the useful notation 
$$
y_\G(z)=\max_\ca (\max_{\g \in \G}(\Im( \sa^{-1}\g z)))
$$
which measures how close $z \in \H$ is to a cusp. If $\psi$ (or $|\psi|$) 
is smooth with weight $0$ then it is more convenient to write
$$
\psi(z) \ll y_\G(z)^A,
$$
for example, instead of $\psi(\sa z) \ll y^A$ for each cusp $\ca$ as $y 
\rightarrow \infty$. 
For questions involving the convergence of inner products it is the 
growth in these cuspidal zones that is important.
For weight $0$ second-order forms it makes more sense to consider their 
growth only inside the fundamental domain $\F$. We use the notation
$$
y_{\F}(z)=\max_\ca (\Im( \sa^{-1} z))
$$
for $z \in \F$.

\vskip .10in Next, we recall from Section 2 the non-holomorphic, weight $k$
Poincar\'e series
$$
U_{\ca m}(z,s,k)=\sum_{\g \in \G_\ca \backslash \G}   \Im(\sa^{-1} \g 
z)^s e(m \sa^{-1} \g z) \ee(\sa^{-1} \g, z)^{-k} \tag 5.1
$$
where $\ee(\g, z)$ is the `non-holomorphic' weight factor
$j(\g,z)/|j(\g,z)|$. These series were first studied by Selberg.
For simplicity put $U_{\ca m}(z,s)=U_{\ca m}(z,s,0)$.
The main facts we need for this function are that it converges to an analytic function of $s$ for  $\Re(s)>1$ and has a continuation to a neighborhood of $s=1$. These results are contained in Propositions A, B and C below.

\vskip .10in Now, for $f(z)$ in $S_2(\G)$ we shall require its derivatives and 
antiderivatives. Define $I_{\ca n}$ to be the $n$th antiderivative of 
$f$. Precisely, for $n \geqslant 1$, we set
$$
I_{\ca n}(z_n)=\int_{i\ci}^{z_n} \cdots 
\int_{i\ci}^{z_2}\int_{i\ci}^{z_1} f_\ca(z_0) dz_0 dz_1 \cdots dz_{n-1}
$$
for $f_\ca(z)=f(\sa z)/j(\sa,z)^2$. Thus $\frac d{dz} I_{\ca n}(z)=I_{\ca 
(n-1)}(z)$ and use this to extend the definition of $I_{\ca n}$ to all 
$n$ in $\Z$. In particular it may be checked that
$$
\align
I_{\ca 1}(\sa^{-1}z)&=F_{\ca}(z):=\int_{\ca}^z f(w) \,dw,\\
I_{\ca 0}(\sa^{-1}z)&=f(z)j(\sa^{-1},z)^2.
\endalign
$$
There are two interesting  families of series:
$$
\align
Q_{\ca m}(z,s,n;f)&=\sum_{\g \in \G_\ca \backslash \G}  I_{\ca 
n}(\sa^{-1} \g z) \Im(\sa^{-1} \g z)^s e(m \sa^{-1} \g z),\\
Q_{\ca m}(z,s,n;\overline{f})&=\sum_{\g \in \G_\ca \backslash \G}  
\overline{I_{\ca n}(\sa^{-1} \g z)} \Im(\sa^{-1} \g z)^s e(m \sa^{-1} \g 
z).\tag 5.2
\endalign
$$
For $m=0$ the theory of the series $Q_{\ca m}(z,s,n;f)$ for all $n$ and 
$f$ with weight $k\geqslant 2$ is given in \cite{CO}. 
In \cite{JO} the series $Q_{\ca m}(z,s,1;f)$ with $m \neq 0$ is used to 
find the analogue of the Kronecker limit formula for second-order, 
non-holomorphic Eisenstein series.
Here we will only need results about $Q_{\ca m}(z,s,n;\overline{f})$ and 
only for $n\leqslant 1$.

\vskip .10in The next proposition provides the basic convergence results and growth estimates for $U_{\ca m}(z,s,k)$, 
$Q_{\ca m}(z,s,1;\overline{f})$, $Q'_{\ca m}(z,s,1;\overline{f})$ which is the termwise derivative of $(5.2)$ with respect to $z$, and the weight 2 function
$$
G_{\ca m}(z,s;\overline{F})=\sum_{\g \in \G_\ca \backslash \G} 
\frac{\overline{F_\ca(\g z)}}{j(\sa^{-1} \g, z)^2}  \Im(\sa^{-1} \g z)^s e(m 
\sa^{-1} \g z).
$$

\proclaim{Proposition A} For $k \in 2\Z$ and $\sigma=\Re (s) >1$ the series
$U_{\ca m}(z,s,k)$, $Q_{\ca m}(z,s,1;\overline{f})$, $Q'_{\ca m}(z,s,1;\overline{f})$ and $G_{\ca m}(z,s-1;\overline{F})$ converge absolutely and uniformly on compact sets to analytic functions of $s$. For these $s$ we
have
$$
\align
U_{\ca 0}(z,s,k) & \ll   y_\G(z)^{\sigma},\tag i\\
U_{\ca m}(z,s,k) & \ll  1, \ \ \ m>0,\tag ii\\
Q_{\ca m}(z,s,1;\overline{f}) & \ll   y_\G(z)^{1/2-\sigma/2}, \ \ \ m \geqslant 0,\tag iii\\
yQ'_{\ca m}(z,s,1;\overline{f}) & \ll   (|m|+1) y_\G(z)^{1/2-\sigma/2}, \ \ \ m \geqslant 0,\tag iv\\
yG_{\ca m}(z,s-1;\overline{F}) & \ll  y_\G(z)^{1/2-\sigma/2}, \ \ \ m \geqslant 0 \tag v
\endalign
$$
where the implied constants depends on $s$, $k$, $f$ and $\G$ but not $m$.
\endproclaim

The proof of this proposition is not difficult and essentially amounts to comparing these series with the standard Eisenstein series $E_{\ca}(z,s)$. We give the details in Section 9. The next propositions show the analytic continuation of the  series $U_{\ca m}$, $Q_{\ca m}$ and $Q'_{\ca m}$ to a neighborhood of $s=1$. This is  essential for our construction of weight two second-order forms. Again, to keep the flow of ideas intact, we relegate the proofs to Sections 10 and 11.

\vskip .10in First we choose, once and for all, a constant $\delta_\G$ depending on $\G$ with $0<\delta_\G <1/2$. It is chosen so that poles appearing from the discrete spectrum have real part less than $1-\delta_\G$. See the discussion at the beginning of Section 8 for a complete explanation.

\proclaim{Proposition B} For $k \in 2\Z$  the (Eisenstein) series
$U_{\ca 0}(z,s,k)$  has a meromorphic
continuation to all $s$ with $\Re (s) >1-\delta_\G$. We
have
$$
U_{\ca 0}(z,s,k) \ll   y_\G(z)^{\sigma}.
$$
for these $s$ with the implied constant depending on $s$, $k$ and $\G$.
The only possible pole in this region appears at $s=1$
when $k=0$. It is a simple pole with residue $1/V$.
\endproclaim

\proclaim{Proposition C} For $k \in 2\Z$ and $m > 0$ the Poincar\'e series
$U_{\ca m}(z,s,k)$  has a
continuation to an analytic function for all $s$ with $\Re (s) >1-\delta_\G$. We
have
$$
U_{\ca m}(z,s,k) \ll   y_\G(z)^{1/2}
$$
for these $s$ with the implied constant depending on $s$, $m$, $k$ and $\G$.
\endproclaim

This proposition should be standard and the proof of the first part can 
be found in [S]. However, since we have not found a reference for the 
growth estimate, we provide a proof in Section 10 based on the methods of [JO] Section 
8.

\proclaim{Proposition D} For $m \geqslant 0$, both series $(s-1)Q_{\ca 
m}(z,s,1;\overline{f})$ and $Q'_{\ca m}(z,s,1;\overline{f})$ have  
continuations to analytic functions of $s$ with 
$\Re (s) >1-\delta_\G$. For these $s$ values
$$
(s-1)Q_{\ca m}(z,s,1;\overline{f}), \ \ \ yQ'_{\ca m}(z,s,1;\overline{f}) \ll  
y_\G(z)^{1/2}.
$$
The implied constant depends on $s$, $m$, $f$ and $\G$. Also $Q_{\ca 
m}(z,s,1;\overline{f})$ has a simple pole at $s=1$ with residue 
$2i\overline{\s{f}{P_{\ca m}(\cdot)_2}}$.
\endproclaim

The $P_{\ca m}(z)_2$ appearing in the above residue is
$$
P_{\ca m}(z)_2=y^{-1}U_{\ca m}(z,1,2), \tag 5.3
$$
which is well defined thanks to Propositions B and C.
These series are holomorphic for $m>0$ and span $S_2(\G)$, see 
\cite{JO} Theorem 3.2 for example. For $m=0$ they are no longer holomorphic. One 
way to see this is to note that
$$
P_{\ca 0}(z)_2=y^{-1}U_{\ca 0}(z,1,2)=y^{-1}\lim_{s \rightarrow 1} R_0 
E_{\ca}(z,s)= 2i \lim_{s \rightarrow 1} \frac{d}{dz} E_{\ca}(z,s)
$$
which by (4.7) satisfies
$$
j(\sb,z)^{-2}P_{\ca 0}(\sb z)_2=\delta_{\ca \cb}-\frac{1}{y V}+O(e^{-2\pi y}) \tag 5.4
$$
as $y\rightarrow \ci$ and
$$
y^2\frac d{d\overline{z}}P_{\ca 0}(z)_2=\frac{i}{2 V}. \tag 5.5
$$
These functions were also used in \cite{GO} where it is shown that for any two distinct cusps $\ca$ 
and $\cb$ the differences
$P_{\ca 0}(z)_2-P_{\cb 0}(z)_2$ span $E_2(\G)$, the space of Eisenstein series (4.2). The operator $R_0=2iy d/dz$ is the weight raising operator as discussed in Section 8.

\vskip 3mm
\flushpar
{\bf Remark.} In reality all the series we are considering here $U_{\ca m}(z,s,k)$, $Q_{\ca m}(z,s,1;\overline{f})$, $yQ'_{\ca m}(z,s,1;\overline{f})$ and 
$yG_{\ca m}(z,s-1;\overline{F})$ have meromorphic continuations to all $s$ in $\C$ and indeed it is a relatively straightforward exercise to extend the proofs of Propositions C and D to the entire $s$ plane. The reason we restrict our attention to $\Re(s)>1-\delta_\G$ is to avoid having to include details concerning poles coming from the exceptional spectrum and on the line $\Re(s)=1/2$. For our purposes we only require continuation to a neighborhood of $s=1$. See \cite{CO} for details of how these techniques extend to all of $\C$

\vskip .10in We are now ready to turn to the dimension of
$S_2^2(\G)$. Because of the decomposition (3.7), in 
order to compute it, we need to investigate the space $S_2^2(\G)^-$. In view of the remarks following Theorem 4.3, we would like to construct the 
series
$$
\text{``}P_{\ca m}(z,L)_2\text{''}=\sum_{\g \in \G_\ca \backslash \G} 
L(\gamma) j(\sa^{-1} \g, z)^{-2} e(m \sa^{-1} \g z)
$$
with $L(\g)=\overline{\s{\g}{f}},$  $f \in S_2(\G).$ Unfortunately this 
series is not absolutely 
convergent. However, if it were, it would formally equal 
$$
\multline
\sum_{\g \in \G_\ca \backslash \G} \overline{(F_\ca(\g z)-F_\ca(z))} 
j(\sa^{-1} \g, z)^{-2} \Im(\sa^{-1} \g z)^s e(m \sa^{-1} \g z)\\
 =G_{\ca m}(z,s;\overline{F})-\overline{F_\ca(z)} y^{-1}U_{\ca m}(z,s+1,2)
\endmultline
$$
with $s=0$. This motivates us to study 
$$
\align
Z_{\ca m}(z,s;f)& = \sum_{\g \in \G_\ca \backslash \G} \overline{\s{\g}{f}} 
\frac{\Im(\sa^{-1} \g z)^s}{j(\sa^{-1} \g, z)^{2}}  e(m \sa^{-1} \g z) \tag 5.6\\
&=G_{\ca m}(z,s;\overline{F})-\overline{F_\ca(z)} 
y^{-1}U_{\ca m}(z,s+1,2) \tag 5.7
\endalign
$$
and show that it has an analytic continuation to $s=0$. An analysis 
will show that it does and has the correct growth at cusps but is not 
always analytic. Nonetheless, this construction provides us with all  
remaining elements of $S_2^2(\G)$ and $M_2^2(\G)$. The next proposition develops the required properties of $Z_{\ca m}$.

\proclaim{Proposition 5.1} 
Let $m$ be an integer $\geqslant 0$, $f \in S_2(\G)$ and $\ca$ a cusp. Then $Z_{\ca m}(z,s;f)$ as defined in $(5.7)$  admits an analytic continuation to  $\Re(s)> -\delta_\G$. Also 
$$
\align 
Z_{\ca m}(\g z,0;f)j(\g,z)^{-2}
&=Z_{\ca m}(z,0;f)-\overline{\s{\g}{f}} P_{\ca m}(z)_2, \tag i\\
yZ_{\ca m}(z,0;f) &\ll y_{\F}(z)^{1/2}, \tag ii\\
\frac d{d\overline{z}}Z_{\ca 
m}(z,0;f)& =-y^{-2}\overline{\s{f}{P_{\ca m}(\cdot)_2}} + \delta_{m,0} \overline{F_\ca (z)}/(2iy^2 V), \tag iii
\endalign
$$
where (i) is true for all $\g$ in $\G$ and the implied constant in (ii) is independent of $z$.
\p We have already seen in Proposition A part (v) that $G_{\ca m}(z,s;\overline{F})$ is absolutely convergent for 
$\sigma=\Re (s)>0$ and satisfies
$$
yG_{\ca m}(z,s;\overline{F}) \ll y_{\G}(z)^{-\sigma/2} \tag 5.8
$$
for these $s$ and an implied constant independent of $z$. Now use the relation
$$
\frac d{dz} Q_{\ca m}(z,s,1;\overline{f})=\frac{-is}{2}G_{\ca 
m}(z,s-1;\overline{F})
+2\pi im G_{\ca m}(z,s;\overline{F})
$$
to see that
$$
G_{\ca m}(z,s;\overline{F})=\frac{4\pi m}{s+1}G_{\ca m}(z,s+1;\overline{F})
+\frac{2i}{s+1} Q'_{\ca m}(z,s+1,1;\overline{f}). \tag 5.9
$$
With Proposition D this gives the analytic continuation of $G_{\ca 
m}(z,s;\overline{F})$ to $\Re (s) >-\delta_\G$. Together with 
Propositions B and C we see that both terms on the right of (5.7) have analytic continuations to $\Re (s) >-\delta_\G$, proving the first statement of the proposition.

\vskip .10in Equations (5.7) and (5.9) now imply that
$$
yZ_{\ca m}(z,0;f) = 4\pi m yG_{\ca m}(z,1;\overline{F})
+2i yQ'_{\ca m}(z,1,1;\overline{f})-\overline{F_\ca(z)} 
U_{\ca m}(z,1,2).
$$
With Proposition A part (v) and Proposition D we have
$$
yG_{\ca m}(z,1;\overline{F}) \ll y_\G(z)^{-1/2}, \ \ 
 yQ'_{\ca m}(z,1,1;\overline{f}) \ll y_\G(z)^{1/2}.
$$
To bound $\overline{F_\ca(z)} U_{\ca m}(z,1,2)$ we first note that for $m \neq 0$ we have $U_{\ca m}(z,1,2) \ll 
y_\G(z)^{1/2}$ by Proposition C and for $m= 0$ we have $U_{\ca m}(\sb z,1,2) \ll 
\delta_{\ca \cb}y +1$ as $y \to \infty$ by (10.3). Also it is easy to show that, for $z$ in $\F$, 
$$
F_\ca(\sa z) \ll e^{-2\pi y}, \ \ F_\ca(\sb z) \ll 1, \ca \neq \cb
$$
as $y \to \infty$ (see (9.3) and (9.4)).
Part (ii) of the proposition now follows.
It is also easy to check that, for $\Re(s)$ large,
$$
Z_{\ca m}(\g z,s;f)j(\g,z)^{-2}=Z_{\ca m}(z,s;f)-\overline{\s{\g}{f}} 
y^{-1}U_{\ca m}(z,s+1,2)
$$
thus deducing (i) by analytic continuation.
A lengthy but routine calculation yields, again for $\Re(s)$ large,
$$
\frac d{d\overline{z}}Z_{\ca m}(z,s;f)
=\frac{is}{2y^2} \left(Q_{\ca 
m}(z,s+1,1;\overline{f})-\overline{F_{\ca}(z)}U_{\ca m}(z,s+1)\right).
$$
Combined with Propositions B, C and D we then find that (iii) holds.
This completes all the parts of Proposition $5.1$. \bs

Before we continue we define a new function, $Z_{\lambda, \mu}$, that is more convenient to work with than $Z_{\ca m}$.
As previously noted, for different values of $m>0$, $P_{\ca m}(\cdot)_2$
spans $S_2(\G)$. So for any pair $(\lambda,\mu)$ of elements in $S_2(\G)$
we may find a linear combination $Z_{\lambda, \mu}$ of
$Z_{\ca m}(z,0;\lambda)$ with  $m>0$ such that, for all $\g$ in $\G$
$$
\align
Z_{\lambda,\mu}|_2 (\g -1) &= \overline{\s{\g}{\lambda}}\mu, \tag 5.10\\
y^2 \frac d{d\overline{z}}Z_{\lambda,\mu}&=\overline{\s{\lambda}{\mu}}, \tag 5.11\\
y Z_{\lambda,\mu} &\ll y_{\F}(z)^{1/2}.\tag 5.12
\endalign
$$

\vskip .10in 
To prove that the dimension of $S_2^2(\G)$ is one less than we would initially expect, we will also need the next result.

\proclaim{Proposition 5.2} For any $\mu \ne 0$ in $S_2(\G)$ it is 
impossible for an element $f$ of $S_2^2(\G)$ to satisfy 
$$
f|_2(\g-1)=\overline{\s{\g}{\mu}}\mu
$$
for all $\g$ in $\G$.
\p
The main idea is to use
$
P_{\ca 0}(z)_2 \in \tilde M_2(\G)
$
and consider the sum $Z_{\mu,\mu}+2iV\overline{\s{\mu}{\mu}} P_{\ca 
0}(\cdot)_2)$. Combine (5.5) and (5.11) to obtain
$$
\frac d{d\overline{z}}(Z_{\mu,\mu}(z)+2iV\overline{\s{\mu}{\mu}} P_{\ca 
0}(z)_2)=0.$$
Looking at each cusp we see by (5.4) and (5.12) that
$$
j(\sb,z)^{-2}(Z_{\mu,\mu}(\sb z)+2iV\overline{\s{\mu}{\mu}} P_{\ca 0}(\sb 
z)_2)\ll y^{-1/2}+\delta_{\ca \cb}+\frac{1}{y V} 
$$
as $y\rightarrow \ci$ so the Fourier expansion of $Z_{\mu,\mu}(z)+2iV\overline{\s{\mu}{\mu}} P_{\ca 
0}(z)_2$ has only non negative terms and only a constant term at the cusp $\ca$. Now suppose $f$, as described in the proposition, 
does exist. Then
$$
\lambda=f-Z_{\mu, \mu}-2iV\overline{\s{\mu}{\mu}} P_{\ca 0}(\cdot)_2
$$
is an element of $M_2(\G)$ that has exponential decay at every cusp 
except possibly one by the previous argument. But $P_{\ca 
0}(\cdot)_2-P_{\cb 0}(\cdot)_2$ span $E_2(\G)$ and 
hence every element of 
$E_2(\G)$ must not have exponential decay at at least two 
distinct cusps.
It follows that $\lambda$ must be in $S_2(\G)$. This yields a 
contradiction on writing
$$
 \lambda(z)-f(z)+Z_{\mu,\mu}(z)=-2iV\overline{\s{\mu}{\mu}} P_{\ca 0}(z)_2
$$
and examining the size of both sides as $z \rightarrow \ca$. \bs

\proclaim{Theorem 5.3} Let $\{f_1, \dots, f_g\}$ be an orthonormal  basis of $S_2(\G)$. 
Then the set
$$
A=
\{Z_{f_i, f_j}\}_{i \ne j} \cup \{Z_{f_i, f_i}-Z_{f_1, f_1}\}_{i=2}^g
\cup \{f_i \int f_j\}_{i, j} \cup \{f_i\}_i \tag 5.13
$$ is a basis of
$S_2^2(\Gamma)$. 
\p
We first observe that, with (5.10), (5.11) and (5.12), all functions  belong to $S_2^2(\G).$
Next we show that $A$ generates $S_2^2(\G).$ Equations (3.3) and (3.5) imply that if $f
 \in S_2^2(\G)$ then there are $\lambda_{ij}, \mu_{ij} \in \Bbb C$ such
that 
$$ f|_2(\g-1)=\sum_{i,
j}(\lambda_{ij} \s{\g}{f_i} +\mu_{ij}\overline{\s{\g}{f_i}})f_j 
$$
for all $\g$ in $\G$.
Set
$$
f^*=f-\sum_{i, j}\lambda_{ij} f_j \int f_i -\sum_{i, j}
\mu_{ij} Z_{f_i, f_j}.
$$
Thanks to  (5.10), 
$f^*|_2(\g -1) =0$ for all $\g$ in $\G$. Writing
$$
f^*=f-\sum_{i, j}\lambda_{ij} f_j \int f_i -\sum_{i \neq j}
\mu_{ij} Z_{f_i, f_j}-\sum_{i}
\mu_{ii} \left( Z_{f_i, f_i}-Z_{f_1, f_1}\right) -\left(\sum_{i}
\mu_{ii}\right) Z_{f_1, f_1}
$$
makes it clear by (5.11) that
$$
y^2 \frac d{d\overline{z}}f^*= -\left(\sum_{i}
\mu_{ii}\right) \overline{\s{f_1}{f_1}}.
$$
Hence
$f^*+(\sum_{i}
\mu_{ii}) Z_{f_1,f_1}$ is holomorphic.

\vskip .10in
With Lemma 4.1, $f_i\int_{z_0}^z \lambda_i(w)\,dw$ satisfies condition
{\bf C}. Also, with $(5.12)$,
$$
j(\sb,z)^{-2}Z_{f_i,f_j}(\sb z)\ll \Im(\sb z)^{-1}
|j(\sb,z)|^{-2}y^{1/2}
=y^{-1/2}
$$
as $y\rightarrow \ci$. Therefore $f^*+(\sum_{i}
\mu_{ii}) Z_{f_1,f_1}$ must have a
Fourier expansion at each cusp with only positive terms. In other words
$$
j(\sb,z)^{-2}\left(f^*(\sb z)+(\sum_{i}
\mu_{ii}) Z_{f_1,f_1}(\sb z)\right)=\sum_{n=1}^\ci d_\cb(n)e^{2\pi in z}
$$
for constants $d_\cb(n)$ and hence $f^*+(\sum_{i}
\mu_{ii}) Z_{f_1,f_1}$ satisfies condition {\bf C}.

\vskip .10in
Further, if $\sum_i \mu_{ii} \ne 0$ then for all $\g \in \G$
$$\left.\frac{f^*+(\sum_{i}
\mu_{ii}) Z_{f_1,f_1}}{\sum_{i}
\mu_{ii}}\right|_k(\gamma-1)=Z_{f_1, f_1}|_k(\g-1)=\overline{\s{\g}{f_1}} f_1$$
which contradicts Proposition 5.2. Consequently $\sum_{i}
\mu_{ii}=0$, $f^* \in S_2(\Gamma)$ and 
$A$ generates $S_2^2(\G)$.

\vskip .10in
Finally, we verify that $A$ is linearly independent. Suppose that for some
$k_{ij}, l_{i j}, n_i, m_i \in \Bbb C$, we have
$$
\sum_{i \ne j} k_{i j} Z_{f_i, f_j}+ \sum_{i, j} l_{i j} f_j \int f_i
+\sum_{i \ne 1} n_i (Z_{f_i, f_i}-Z_{f_1, f_1})+\sum_i m_i f_i=0.$$
With (5.10), if we let $\g-1$, ($\g \in \G$) act on both sides, we 
obtain:
 $$\sum_{i \ne j} k_{i j} \overline{\s{\g}{f_i}} f_j+
\sum_{i, j} l_{i j} \s{\g}{f_i}  f_j+\sum_{j \neq 0} n_j (\overline{\s{\g}{f_j}}f_j-\overline{\s{\g}{f_1}}f_1)=0.$$ 
Take the inner product of both sides with $f_j$, $j \neq 1$ to see that
$$
\sum_{i \ne j} k_{i j} \overline{\s{\g}{f_i}} +
\sum_{i} l_{i j} \s{\g}{f_i}  + n_j \overline{\s{\g}{f_j}}=0.
$$
Therefore $k_{ij}$, $l_{ij}$ and
$n_j$ are all $0$ for $j \neq 1$ by the usual
Eichler-Shimura isomorphism. Similarly for $j=1$. Hence the constants $m_i$ must vanish too. \bs

\proclaim{Corollary 5.4}
$$
\dim S^2_2(\G)=\cases 0 \text{ if }\dim S_2(\G)=0, \\ (2g+1)\dim S_2(\G) -1 \text{ otherwise} \endcases.
$$

\endproclaim

\flushpar
{\bf Remark 5.5.} As a byproduct of the proof of Proposition 5.2, 
we have also shown that the interesting form 
$$
Z_{\mu,\mu}+2iV\overline{\s{\mu}{\mu}} P_{\ca 0}(\cdot )_2
$$
exists and is in $R_2^2(\G)$ but not $S_2^2(\G)$.

$$\text{\bf 6. Calculating $\dim M^2_k(\G)$ for $k \leqslant 2$}$$

First we work on the case $k=2$.
For any $f \in S_2(\G)$ and cusps $\ca,$ $\cb$ put 
$$
\Cal Z_{\ca \cb}(z;f) = Z_{\ca 
0}(z,0;f)-Z_{\cb 0}(z,0;f)+P_{\ca 0}(z)_2 \overline{\int_\ca^\cb f(w)\, dw}.
$$
 With Proposition 5.1, part (iii) and (5.5), we have
$$
-y^2 \frac{d}{d \overline{z}} \Cal Z_{\ca \cb}(z;f)= 
\overline{\s{f}{P_{\ca 0}(\cdot)_2 - P_{\cb 0}(\cdot)_2}}.
$$
It is routine to verify by unfolding that $\s{f}{\Im(\cdot)^{-1}U_{\ca 
0}(\cdot, s,2)} = 0$ for $\Re(s)$ large. It follows by analytic 
continuation that $\s{f}{P_{\ca 0}(\cdot)_2} = 0$ and hence that
$$
\frac{d}{d \overline{z}} \Cal Z_{\ca \cb}(z;f)= 0.
$$
So $\Cal Z_{\ca \cb}(z;f)$ satisfies {\bf H} and, with (ii) of 
Proposition 5.1 and (5.4), it satisfies condition {\bf N}. By part (i) of the same 
proposition
$$
\Cal Z_{\ca \cb}(\g z;f)j(\g, z)^{-2} = \Cal Z_{\ca 
\cb}(z;f)-\overline{\s{\g}{f}} \left(
P_{\ca 0}(z)_2 - P_{\cb 0}(z)_2 \right)
$$
and conditions {\bf A(}$M_k(\G)${\bf )} and {\bf P} hold. We have shown that $\Cal 
Z_{\ca \cb}(z;f) \in M_2^2(\G)$ for all cusps $\ca$, $\cb$ and cusp forms 
$f$. As we already noted, the differences $P_{\ca 0}(z)_2 - P_{\cb 
0}(z)_2$ span $E_2(\G)$ (and are orthogonal to $S_2(\G)$). 
So for any $f\in S_2(\G)$ and $h\in E_2(\G)$ we may find a linear combination, $Z_{f,h}$, of $\Cal Z_{\ca \cb}(z;f)$ for different $\ca$, $\cb$ that is holomorphic and satisfies
$$
Z_{f,h}|_2 (\g -1) = \overline{\s{\g}{f}}h. 
$$
The next theorem is easy to check using Remark 5.5 and the corollary follows quickly.
 
\proclaim{Theorem 6.1} Let 
$\{f_1, \dots, f_g\}$ be an orthonormal  basis of $S_2(\G)$ and $\{h_1, \dots ,h_{p-1}\}$ a basis for 
$E_2(\G)$. If $\frak a$ is a fixed cusp, then
the set
$$
\{Z_{f_i, f_j}\}_{i \ne j} \cup \{Z_{f_i, f_i}+2iVP_{\frak a 0}(\cdot)_2\}_i
\cup \{Z_{f_i,h_j}\}_{i,j} \cup \{f_i \int f_j\}_{i, j} \cup 
\{h_j \int f_i\}_{i,j} \cup \{f_i\}_i  
\cup \{h_j \}_{j}
$$
is a basis for $M_2^2(\G)$.
\endproclaim
 
\proclaim{Corollary 6.2} We have
$$
\dim M^2_2(\G)= (2g+1)\dim M_2(\G).
$$
\endproclaim

Next we look at $M^2_0(\G)$. 

\proclaim{Theorem 6.3} With $\{f_1, \dots, f_g\}$ a basis of $S_2(\G)$,  a basis for $M^2_0(\G)$ is 
$$
\{1\} \cup \{\int f_i\}_i \tag 6.1
$$
and $\dim M^2_0(\G)=g+1$.
\p
Let $f \in M_0^2(\G).$ Since $M_0(\G)$ consists of the constant functions, for each $\g \in \G,$
$f|_0(\g-1)=c_{\g} \in \C$. However $f$
belongs to $M_0^2(\G)$, so the map $\g \to
c_{\g}$ is a parabolic $1$-cocycle in terms of the trivial action of $\G$
on $\Bbb C.$  Hence, by the Eichler-Shimura isomorphism, there are 
$g, h \in S_2(\G)$ such that $c_{\g}=\s{\g}{g}+\overline{\s{\g}{h}}$. This is
equivalent to 
$$
\left. \left(f(z)-\int_{z_0}^z g(w)dw\right)\right|_0(\g-1)=\overline{\s{\g}{h}}. \tag 6.2$$
Set $F(z)=f(z)-\int_{z_0}^z g(w)dw .$ An
easy computation shows that $F'=f'-g$ has weight 2.
Using the Fourier expansion of $f$ at any cusp we also deduce that
$F'$ satisfies condition \bf C \rm . Hence $F' \in S_2(\G)$ and $$(F|_0(\g-1))(z)
=F(\g z)-F(z)=\int_{z}^{\g z} F'(w) dw=\s{\g}{F'}.$$
Because of the Eichler-Shimura isomorphism, this together with (6.2) implies
that $F'=h \equiv 0.$ Hence, since $f$ is holomorphic, $f(z)=\int_{z_0}^z g(w)dw+c,$ where $c$
is a constant. This shows that (6.1) spans $M^2_0(\G)$. Linear independence follows as in the proof of Theorem 5.3.
 \bs

It is easy to check that when the first-order space is 0 that the corresponding second-order space must be 0. For  even $-k \leqslant -2$ we have $M_{-k}(\G)=0$. Hence any $f \in M^2_{-k}(\G)$ satisfies $f|_{-k}(\g -1)=0$ for all $\g$ in $\G$. Therefore $f \in M_{-k}(\G)$ and $f \equiv 0$. The same argument shows that $S^2_{-k}(\G)=0$ for all even $-k \leqslant 0$. 

\vskip .10in
All parts of Theorems 2.1 and 2.2 are now complete.

$$\text{\bf 7. An Eichler-Shimura-like isomorphism.}$$

For $k \geqslant 2$ let $P_{k-2}$ 
denote the space of polynomials of degree at most
$k-2$ with coefficients in $\Bbb C.$ With $F$ in $M^2_k(\G)$ define a map 
$\phi:\Gamma \to P_{k-2}$ by the formula
$$
\phi(\gamma)=\int_{i}^{\gamma^{-1} i}
F(z)(z-X)^{k-2} dz
$$
for all $\gamma \in \Gamma$ where $X$ is the polynomial variable and 
integration takes place on
a geodesic in the upper-half plane. This is the exact analogue of the
period polynomial map of Eichler cohomology and, as in the classical case,
it has a close relation to the values of $L(s, F)$ at $s=1, \dots,
k-1.$ Specifically, if $F$ is cuspidal, $i$ can be replaced by $i \infty$
in the definition of $\phi$ and, then the polynomial coefficients are linear
combinations of values of the additively twisted $L$-function of $F$ which is
$$
\sum_{n=1}^{\infty}\frac{a_n e^{2\pi i mn}}{n^s}
\text{ \ \ for \ \ }F(z)=\sum_{n=1}^{\infty} a_n e^{2\pi i n z}.
$$

\vskip .10in
As usual, for every $\C[\Gamma]$-module $M$ we let $d$ denote the coboundary 
operator on the group of $i$-cochains 
$C^i(\Gamma, M):=\{f: \Gamma^i \to M \}$. For example, for every $\psi: 
\Gamma \to M,$ 
$$
(d\psi)(\gamma_1, \gamma_{2})=\psi(\gamma_2).\gamma_1-
\psi(\gamma_2\gamma_1)+\psi(\gamma_1)
$$ where $.$ denotes the action of 
$\Gamma$ on $M.$ Also, for $\psi: \Gamma^2 \to M,$
$$(d\psi)(\gamma_1, \gamma_{2}, \gamma_3)=\psi(\gamma_2, 
\gamma_3).\gamma_1- \psi(\gamma_2\gamma_1, \gamma_3)+\psi(\gamma_1, 
\gamma_3 \gamma_2)-\psi(\gamma_1, \gamma_2).$$ We 
write $Z^i(\Gamma, M),$  $B^i(\Gamma, M)$ and $H^i(\Gamma, M)$ for the 
groups of $i$-cocycles, $i$-coboundaries and $i$-cohomology 
classes respectively. See, for example \cite{Sh} p. 223 for more details.

\vskip .10in
Further, we set $C^q_{par}(\Gamma, M)=C^q(\Gamma, M)$ for $q \ne 1$ and 
we we define the group of parabolic $1$-cochains 
$$
C^1_{par}(\Gamma, M)=\{f \in C^1(\Gamma, M) \ | \ f(\pi) \in M.(\pi-1) 
\text{ for all parabolic }\pi \text{ in }\G \}.
$$
In this way we obtain a cochain complex in terms of the usual coboundary 
operator. Therefore, we have a (``parabolic") cohomology which we can 
describe explicitly as follows: \newline
Set $Z^1_{par}(\Gamma, M)=Z^1(\Gamma, M) \cap C^1_{par}(\Gamma, M)$
and $B^2_{par}(\Gamma, M)=d(C^1_{par}(\Gamma, M)).$ With this notation 
define $H^1_{par}(\Gamma, M)=Z^1_{par}(\Gamma, M)/B^1(\Gamma, M)$
and $H^2_{par}(\Gamma, M)=Z^2(\Gamma, M)/B^2_{par}(\Gamma, M).$
In
particular, for $M=\Bbb C$ with the trivial action of $\Gamma$,
$H^1_{par}(\Gamma, \Bbb C)$ (and $Z^1_{par}(\Gamma, \Bbb C)$) is
isomorphic to the group Hom$_0(\Gamma, \Bbb C)$, defined in Section 3.

\vskip .10in
Let now $\Gamma$ act on $P_{k-2}$ via $|_{2-k}$ and on $C^1(\Gamma, 
P_{k-2})$ via the trivial action. Consider the map $\alpha: C^1(\Gamma, 
P_{k-2}) \to C^1(\Gamma, C^1(\Gamma, P_{k-2}))$ defined as follows:
For $\psi: \Gamma \to P_{k-2}$ we let $\alpha(\psi)$ be a map from 
$\Gamma$ to $C^1(\Gamma, P_{k-2})$ such that $\alpha(\psi)(\gamma)$, for 
$\gamma$ in  $\Gamma$, is defined by the formula
$$
\alpha(\psi)(\gamma)(\delta)=(d\psi)(\gamma, 
\delta)|_{2-k}\gamma^{-1}
$$ 
for all $\delta 
\in \Gamma$. 
We then set 
$$\align
Z^1_!(\Gamma, P_{k-2})&=\alpha^{-1}(H^1_{par}(\Gamma, Z^1(\Gamma, 
P_{k-2})))\\
B^1_!(\Gamma, P_{k-2})&=\alpha^{-1}(H^1_{par}(\Gamma, B^1(\Gamma, 
P_{k-2}))).
\endalign
$$
Explicitly, a map $f: \Gamma \to P_{k-2}$ belongs to $Z^1_!(\Gamma, 
P_{k-2})$ if and only if 
\roster
\item"(i)" for each $\gamma \in \Gamma,$ 
$(df)(\gamma, \delta)|_{2-k}\gamma^{-1}$ is a $1$-cocycle as a 
function of $\delta \in \G$, 
\item"(ii)" $(df)(\gamma_2 \gamma_1, 
\delta)|_{2-k}(\gamma_2 \gamma_1)^{-1}=
(df)(\gamma_2, \delta)|_{2-k}\gamma_2^{-1}+
(df)(\gamma_1, \delta)|_{2-k}\gamma_1^{-1}$ for all $\g_1$, $\g_2$, $\delta$ in $\G$ and 
\item"(iii)" $(df)(\pi, \delta) \equiv 0$ for  parabolic $\pi$ and all $\delta$ in $\G$. 
\endroster
In other words $f \in Z^1_!(\Gamma, P_{k-2})$ if and only if 
$$
\align
f(\gamma_3 \gamma_2 \gamma_1)= & \ f(\gamma_3 \gamma_2)|_{2-k} 
\gamma_1 +f(\gamma_2 \gamma_1) +f(\gamma_3 \gamma_1)|_{2-k} 
(\gamma_1^{-1} \gamma_2 \gamma_1)
\\
& -f(\gamma_3)|_{2-k}(\gamma_2 
\gamma_1)-f(\gamma_2)|_{2-k} \gamma_1 - f(\gamma_1)|_{2-k}
(\gamma_1^{-1} \gamma_2 \gamma_1)
\endalign
$$  
and
$$
f(\delta \pi)=f(\delta)|_{2-k}\pi+ f(\pi).
$$
Exactly the same is true for $f \in B^1_!(\Gamma, P_{k-2})$ except that, 
in addition, $(df)(\gamma, \delta)|_{2-k}\gamma^{-1}$ is a $1$-coboundary.
Therefore $B^1_!(\Gamma, P_{k-2})$ can be canonically embedded in 
$Z^1_!(\Gamma, P_{k-2})$
and we may define
$$H^1_{!}(\Gamma, P_{k-2})=\frac{Z^1_!(\Gamma,
P_{k-2})}{B^1_!(\Gamma, P_{k-2})}.$$
We also set $\bar S_k^2(\Gamma)$ (resp. $\bar S_k(\Gamma)$) for the
space of functions whose complex conjugate is in $S_k^2(\Gamma)$
(resp. $S_k(\Gamma)$) and we associate a map $\tilde
\phi$ to $F \in \bar S_k^2(\Gamma)$ by setting
$$\tilde \phi(\gamma)=\int_i^{\gamma^{-1} i}
F(z)(\bar z-X)^{k-2} d \bar z \qquad \text{for all} \, \, \, \gamma \in
\Gamma.$$
With this notation we have

\proclaim{Theorem 7.1}
(i) For $k>2$ the maps defined by $\phi$ and $\tilde \phi$ induce an
isomorphism
$$\frac{M_k^2(\Gamma)}{M_k(\Gamma)} \oplus \frac{\bar S_k^2(\Gamma)}{\bar
S_k(\Gamma)} \cong H^1_{!}(\Gamma, P_{k-2} ).$$
(ii) There is an isomorphism
$$\frac{M_2^2(\Gamma)}{M_2(\Gamma)} \oplus \frac{\bar S_2^2(\Gamma)}{\bar
S_2(\Gamma)} \oplus \Bbb C \cong H^1_{!}(\Gamma, \Bbb C ).$$
\p
(i) We will be using the following formulation of Eichler-Shimura's
isomorphism (see [DI], Section 12.2 for a similar formulation):

\vskip .10in    
\it For each $\phi \in Z^1(\Gamma,
P_{k-2})$ there is a unique pair $(g_1, \bar g_2) \in M_k \times \bar
S_k$ such that
$$
\phi(\gamma)=\int_i^{\gamma^{-1} i}g_1(w)(w-X)^{k-2}dw+
\int_i^{\gamma^{-1} i}\bar g_2(w)(\bar w-X)^{k-2}d\bar w \qquad \text{for all} \, \, \, 
\gamma \in \Gamma.
$$ 
Moreover, if $f$ is the map sending $\phi$ to $(g_1,
\bar g_2)$, then the sequence
$$0\rightarrow B^1(\Gamma, P_{k-2}) \overset{i} \to
\hookrightarrow Z^1(\Gamma, P_{k-2})
\overset {f} \to \rightarrow M_k \oplus \bar S_k \rightarrow 0 \tag7.1$$  
is exact. \rm

\vskip .10in The exact sequence (7.1) can be thought of as a sequence of
$\Gamma$-modules
with $\Gamma$ acting trivially on each of the modules.
This induces an exact sequence of cochain complexes
$$0 \rightarrow C^*_{par}(\Gamma, B^1(\Gamma, P_{k-2})) 
\hookrightarrow C^*_{par}(\Gamma, Z^1(\Gamma, P_{k-2})) 
\rightarrow
C^*_{par}(\Gamma, M_k \oplus \bar S_k) \rightarrow 0$$
and this, in turn, induces a long exact sequence
$$\align
&
 H^1_{par}(\Gamma, B^1(\Gamma, P_{k-2}))
\overset{i^*} \to \rightarrow H^1_{par}(\Gamma, Z^1(\Gamma, P_{k-2}))
\overset{f^*} \to \rightarrow H^1_{par}(\Gamma, M_k \oplus \bar S_k) \\
 &\rightarrow H^2_{par}(\Gamma, B^1(\Gamma, P_{k-2}))
\to H^2_{par}(\Gamma, Z^1(\Gamma, P_{k-2})) \to H^2_{par}(\Gamma,
M_k \oplus \bar S_k) \\
 &\to H^3(\Gamma, B^1(\Gamma,
P_{k-2})) \rightarrow \dots \tag7.2
\endalign
 $$
where $f^*(\psi)(\gamma):=f(\psi(\gamma))$ and $i^*$ is induced by the
injection
$i$ in a similar manner.

\proclaim{Lemma 7.2} $H^j(\Gamma, M)=0$ for every $j \geqslant 2$ and every 
$\Bbb 
C$-vector space $M.$ 
\p 
By Selberg's Lemma (cf. [Ra] or [Se]), there exists a torsion-free subgroup 
$G$ of finite index in $\Gamma$. Since $G$ is itself a Fuchsian group of 
the first kind, it can be
described through the classical generators and relations. Being 
torsion-free
implies that the only relation is 
$$[\gamma_1, \gamma_{g+1}]\dots 
[\gamma_{g}, \gamma_{2g}]\pi_1 \dots \pi_{p}=1,
$$ where $\gamma_j$ (resp. 
$\pi_l$)
are hyperbolic (resp. parabolic) generators of $G$. Therefore, $G$ is 
freely generated by $\gamma_1, \dots \gamma_{2g}$ and $\pi_1, \dots, 
\pi_{p-1}.$
Since the cohomological dimension of a free group is $1,$ an
application of the transfer operator ([B], Proposition 10.1 of III) 
implies that $H^j(\Gamma, M)=0$ for all $j \geqslant 2.$ This completes the proof of the lemma. \bs

We also have (see [Sh], Chapter 8) 
$$H^2_{par}(\Gamma, M) \cong M/M_1$$
where $M_1$ is the subspace of $M$ generated by $M.(\gamma-1)$ for all 
$\gamma$ in $\Gamma.$
Therefore (7.2) implies the exact sequence
$$
\align 
 0 &\to i^*(H^1_{par}(\Gamma,  B^1(\Gamma, P_{k-2})))  
\hookrightarrow H^1_{par}(\Gamma, Z^1(\Gamma, P_{k-2})) 
\overset{f^*} \to \rightarrow H^1_{par}(\Gamma, M_k \oplus \bar S_k) \\
& \rightarrow B^1(\Gamma, P_{k-2}) 
\to Z^1(\Gamma, P_{k-2}) \to M_k \oplus \bar S_k \to 0 .
\endalign 
$$
In particular, this implies that 
$$\dim(i^*(H^1_{par}(\Gamma,  B^1(\Gamma, 
P_{k-2}))))-\dim(H^1_{par}(\Gamma, Z^1(\Gamma, 
P_{k-2})))+\dots-\dim(M_k 
\oplus \bar S_k)=0.$$ Because of (7.1), the last three dimensions cancel 
out 
and we are left with 
$$
\dim(H^1_{par}(\Gamma, Z^1(\Gamma,P_{k-2})))-
\dim(i^*(H^1_{par}(\Gamma,  B^1(\Gamma, P_{k-2}))))=
\dim(H^1_{par}(\Gamma, M_k \oplus \bar S_k)).
$$ 
From this we conclude that $f^*$ induces an isomorphism
$$\frac{H^1_{par}(\Gamma, Z^1(\Gamma, P_{k-2}))}{i^*(H^1_{par}(\Gamma, 
B^1(\Gamma, P_{k-2})))} \cong H^1_{par}(\Gamma, M_k \oplus \bar S_k).$$
Because of the relations (2.1) and 
the vanishing of homomorphisms at elements of
finite order, each element of $H^1_{par}(\Gamma, M_k \oplus \bar S_k)$ is 
uniquely determined by its values at 
the $2g$
hyperbolic generators of $\Gamma$.
Therefore,  
$$
H^1_{par}(\Gamma, M_k \oplus \bar S_k) \cong \bigoplus_{i=1}^{2g} (M_k 
\oplus \bar S_k) \tag 7.3
$$
On the other hand, 
as it was pointed out in [CDO], we have exact
sequences
$$0 \hookrightarrow S_k \to S_k^2 \overset{\psi} \to \rightarrow
\bigoplus_{i=1}^{2g}
S_k$$
$$0 \hookrightarrow M_k \to M_k^2 \overset{\psi} \to \rightarrow
\bigoplus_{i=1}^{2g} M_k \tag7.4$$
where $\psi$ sends a form $F$ to the vector $(F|_k(\gamma_1^{-1}-1), 
\dots, F|_k(\gamma_{2g}^{-1}-1)),$ ($\gamma_i$'s being the hyperbolic
generators of $\G$ we fixed in Section 2).
Theorems 2.1 and 2.2 allow us now to show, by comparison of
dimensions, that, in addition, $\psi$ is onto in both exact sequences
for $k>2$. Therefore, the right-hand side of (7.3) is isomorphic
to $M_k^2 /M_k \oplus \bar S_k^2 / \bar S_k$ via the map
$\psi \times \tilde \psi$, where $\tilde \psi(f)$ is
defined by $\tilde \psi(f):=\overline{\psi(\bar f)}$ (the complex  
conjugation on the right-hand side being understood component-wise).
Therefore,
$$\frac{H^1_{par}(\Gamma, Z^1(\Gamma, P_{k-2}))}{i^*(H^1_{par}(\Gamma, 
B^1(\Gamma, P_{k-2})))}
 \cong H^1_{par}(\Gamma, M_k \oplus \bar S_k) \cong 
\frac{M_k^2}{M_k} \oplus
\frac{\bar S_k^2}{\bar S_k}. \tag7.5$$

Finally, 
\proclaim{Lemma 7.3} The sequence
$$
0 \to Z^1(\Gamma, P_{k-2}) \hookrightarrow Z^1_!(\Gamma, P_{k-2}) 
\overset {\alpha} \to
\rightarrow H^1_{par}(\Gamma, Z^1(\Gamma, P_{k-2})) \rightarrow
0
$$
is exact. 
\p
Since $\alpha(\psi) \equiv 0$ if and only if $d\psi \equiv
0,$ ker$(\alpha)=Z^1(\Gamma, P_{k-2}).$ On the other hand, by the 
definition of $Z^1_!(\Gamma, P_{k-2}),$ we have  
im$(\alpha) \subset H^1_{par}(\Gamma, Z^1(\Gamma, P_{k-2})).$
If $\chi \in H^1_{par}(\Gamma, Z^1(\Gamma, P_{k-2})),$
an easy computation implies that $\chi(\gamma)(\delta)|_{2-k}\gamma$
gives rise to a $2$-cocycle.
From Lemma 7.2, $H^2(\Gamma,
P_{k-2})=0$ and thus, there exists a $\phi:\Gamma \to P_{k-2}$ such that
$(d\phi)(\gamma, \delta)|_{2-k}\gamma^{-1}=\chi(\gamma)( \delta),$ i.e. 
$\chi=\alpha(\phi)$ and $\phi \in Z_!^1.$ This
shows that $H^1_{par}(\Gamma, Z^1(\Gamma, P_{k-2})) \subset$ im$(\alpha).$
\bs

In the same way we can show the exactness of the sequence:
$$0 \to Z^1(\Gamma, P_{k-2}) \hookrightarrow B^1_!(\Gamma, P_{k-2}) 
\overset {\alpha}
\to \rightarrow H^1_{par}(\Gamma, B^1(\Gamma, P_{k-2})) \rightarrow 0.$$
The last two sequences in combination with (7.5) imply the 
isomorphism 
$$\frac{M_k^2}{M_k} \oplus
\frac{\bar S_k^2}{\bar S_k} \cong 
\Big ( \frac{Z^1_!(\Gamma, P_{k-2})}{Z^1(\Gamma, P_{k-2})} \Big )\Big /
\Big ( \frac{B^1_!(\Gamma, P_{k-2})}{Z^1(\Gamma, P_{k-2})} \Big )
\cong 
\frac{Z^1_!(\Gamma, P_{k-2})}{B^1_!(\Gamma, P_{k-2})}.$$

To show that this isomorphism is induced by $\phi \times \tilde \phi$ we 
unravel the
definitions of our maps: $F \in M^2_k$ is first mapped, via $\psi$ to 
$(F|_k(\gamma_i-1))_{i=1}^{2g},$ or, equivalently, to $\gamma \mapsto 
F|_k(\gamma^{-1}-1)$. This is mapped to a $F_1 
\in$ Hom$(\Gamma, Z^1(\Gamma, P_{k-2}))$ such that $f^*(F_1)(\gamma)=
F|_k(\gamma^{-1}-1)$. By  the definition of $f^*$ we have
$f(F_1(\gamma))=F|_k(\gamma^{-1}-1)$ and hence
$$F_1(\gamma)(\delta)=
\int_{i}^{\delta^{-1}i}(F|_k(\gamma^{-1}-1))(w)(w-X)^{k-2}dw.$$
To determine the image of $F_1$ in $H^1_!$
we identify an element $\psi \in Z^1_!$ such that
$(d\psi)(\gamma, \delta)|_{2-k}\gamma^{-1}=
F_1(\gamma)(\delta).$ We verify that the $\phi$ given in the beginning of 
the section is such a map as follows.
The change of variables $w \to \gamma^{-1}  w$
in the integral
$$\phi(\delta \gamma)-\phi(\gamma)=\int_{\gamma^{-1}i}^{\gamma^{-1} 
\delta^{-1} i} F(w)(w-X)^{k-2}dw$$ gives
$$
\int_{i}^{\delta^{-1} i} F(\gamma^{-1} w) (
\gamma^{-1}w-X)^{k-2}d(\gamma^{-1}w)=
\int_{i}^
{\delta^{-1} i} (F|_{k} \gamma^{-1})(w) j(\gamma^{-1}, 
w)^{k-2}(\gamma^{-1}w-X)^{k-2}d w 
$$
where we used $d(\g w)=\frac{dw}{j(\g, w)^2}$. This, in 
combination
with the identity 
$$
(w-\g X)j(\g, X)=(\g^{-1}w-X)j(\g^{-1}, w),
$$
 implies
$$
(d\phi)(\gamma, \delta)=\phi(\delta 
\gamma)-\phi(\delta)|_{2-k}\gamma-\phi(\gamma)=\Big [ \int_{i}^
{\delta^{-1} i} F|_{k} (\gamma^{-1}-1)(w) 
(w-X)^{k-2}d w\Big 
]\Big |_{2-k}\gamma. 
$$
We work in a similar way for $\tilde \phi$ and $\bar S^2_k.$

\vskip .10in
(ii) We work in the same way for $k=2.$ The only difference is at (7.5).
Specifically, because of Theorems 2.1 and 2.2, the exact sequences (7.4)
imply 
$$\frac{M_2^2}{M_2} \oplus \left (
\frac{\bar S_2^2}{\bar S_2} \oplus
\Bbb C \right ) \cong \bigoplus_{i=1}^{2g} (M_2 \oplus \bar S_2).$$
Therefore,
$$\frac{H^1_{par}(\Gamma, Z^1(\Gamma, \Bbb C))}{i^*(H^1_{par}(\Gamma,
B^1(\Gamma, \Bbb C)))}
 \cong H^1_{par}(\Gamma, M_2 \oplus \bar S_2) \cong
\frac{M_2^2}{M_2} \oplus
\frac{\bar S_2^2}{\bar S_2} \oplus \Bbb C.$$
As in the proof of part (i), this implies 
the desired isomorphism and completes the proof of Theorem 7.1.
\bs

\newpage

$$\text{\bf 8. Spectral Theory}$$

Before proceeding to the proofs of Propositions A, B, C and D, 
in Sections 9, 10 and 11, we gather here results we shall require from 
the spectral theory of automorphic forms.

\vskip .10in
We will study the series $U_{\ca m}(z,s,k)$ and $Q_{\ca m}(z,s,1;\overline{f})$ by means of their spectral expansions, which we now recall
(see, for example, \cite{Iw1} and references therein for further
background information and complete proofs). The hyperbolic
Laplacian $\Delta=-4y^2 \ d/dz \ d/d\overline{z}$ operates on $L^2(\G \backslash \H)$ the space of smooth, automorphic, square integrable functions.
Any element $\xi$ of $L^2(\G \backslash \H)$ may be decomposed
into constituent parts from the discrete and continuous spectrum
of $\Delta$. This Roelcke-Selberg decomposition amounts to the identity
$$
\xi(z)=\sum_{j=0}^\infty\langle \xi,\eta_j\rangle
\eta_j(z)+\frac{1}{4\pi }\sum_\cb \int_{-\infty}^{\infty}\langle
\xi ,E_\cb(\cdot,1/2+ir)\rangle E_\cb(z,1/2+ir)\,dr,\tag 8.1
$$
where $\{\eta_j\}$ denotes a complete orthonormal basis of Maass
forms, with corresponding eigenvalues $\lambda_j=s_j(1-s_j)$,
which forms the discrete spectrum.  For notational convenience, we
write $\langle \cdot, \cdot \rangle$ here for the inner product on $\Gamma \backslash \H$ of
weight zero forms (i.e. $\G$-invariant functions).  As always, we
will write $s_j=\sigma_j+it_j$, chosen so that $\sigma_j \geqslant
1/2$ and $t_j \geqslant 0$, and we enumerate the eigenvalues,
counted with multiplicity, by $0=\lambda_0< \lambda_1 \leqslant
\lambda_2 \leqslant \cdots $. Recall that Weyl's law ((11.3) of \cite{I1}) implies
$$
\#\{j| \ |\lambda_j| \leqslant T\} \ll T. \tag 8.2
$$
The decomposition (8.1) is absolutely convergent for each fixed $z$ and
uniform on compact subsets of $\H$, provided $\xi$ and $\Delta
\xi$ are smooth and bounded (see, for example, Theorem 4.7 and
Theorem 7.3 of \cite{I1}). 

\vskip .10in
For each $j$, the Fourier expansion
of $\eta_j$ is
$$
\eta_j(\sa z)=\rho_{\ca j}(0)y^{1-s_j}+\sum_{m\neq 0}
\rho_{\ca j}(m)W_{s_j}(mz).\tag 8.3
$$
For all but finitely many of the $j$ (corresponding to
$\lambda_{j} < 1/4$) we have $\sigma_j=1/2$ and $\rho_{\ca
j}(0)=0$.   
The constant $\delta_\G$ used in Propositions B, C and D, and throughout this paper is chosen so that $1-\delta_\G> \sigma_1 \geqslant 1/2$.

\vskip .10in
We have the bounds
$$
\align
\rho_{\ca j}(0) & \ll 1, \tag 8.4\\
\rho_{\ca j}(m) & \ll \frac{|t_j|}{\sqrt{|m|}}e^{\pi |t_j|/2},  \ \ \ m \neq 0 \tag 8.5
\endalign
$$
where the implied constant depends only on $\G$. The estimate $(8.4)$ is true because, again, there are only finitely many $j$ with $\rho_{\ca j}(0) \neq 0$. The estimate $(8.5)$ follows from the formula of Bruggeman and Kuznetsov, as stated in (9.13) of \cite{Iw1}. See $(8.6)$ of \cite{JO}  for the simple derivation of $(8.5)$. We also have, for any $k \geqslant 0$
and $\sigma=$Re$(s)>1/2-k$,  the bound
$$
W_s(nz)\ll \frac{|s|^{2k}+1}{(|n|y)^{2k-1+\sigma}}|\G(s)|,\tag 8.6
$$
from \cite{JO}, $(8.11)$ (with implied constant depending solely on $\sigma$ and $k$) and Stirling's classical formula
$$
|\G(\sigma+it)| \sim \sqrt{2\pi} |t|^{\sigma-1/2} e^{-\pi|t|/2}
\,\,\,\,\,\text{\rm as $|t|\rightarrow \infty$}. \tag 8.7
$$
Next we bound $\eta_j(z)$. Combine $(8.4)$, $(8.5)$, $(8.6)$, $(8.7)$ and the Fourier expansion $(8.3)$ to obtain
$$
\eta_j(\sa z) \ll y^{1/2} +(|t_j|^{7/2} +1) y^{-3/2} \sum_{m \neq 0} |m|^{-2}
$$
as $y \mapsto \ci$. Therefore
$$
\eta_j(z) \ll y_\G (z)^{1/2} +(|t_j|^{7/2} +1) y_\G (z)^{-3/2} \tag 8.8
$$
for an implied constant depending on $\G$ alone. This is $(8.13)$ of \cite{JO}.

\vskip .10in
Recall the Fourier expansion of $E_{\ca}(z,s)$ in (4.5). The analogues of (8.4) and (8.5) for the continuous spectrum are (we will always assume $T \geqslant 0$ for simplicity)
$$
\align
|\phi_{\ca \cb}(1/2+ir)| & \leqslant 1, \tag 8.9\\
\int_{T}^{T+1} |\phi_{\ca \cb}(m,1/2+ir)|^2\, dr &\ll \frac{T^2}{|m|}e^{\pi T}. \tag 8.10
\endalign
$$
where (8.9) follows from \cite{Iw} (6.28) and (8.10) again from the Bruggeman and Kuznetsov formula as in  \cite{JO} (8.7). Another useful direct bound, as shown in \cite{JO} Lemma 8.4, is 
$$
\phi_{\ca \cb}(m,1/2+ir)  \ll |m|^2, \tag 8.11
$$
for $r$ in $[T,T+1]$ and an implied constant depending on $T$ and $\G$ alone.

\vskip .10in
Let
$C^\infty(\G\backslash\H,k)$ denote the space of smooth functions
$\psi$ on $\H$ that transform as
$$
\psi(\g z)=\ee(\g,z)^k\psi(z)
$$
for $\g$ in $\G$ and $\ee(\g,z)=j(\g,z)/|j(\g,z)|$. For example
$
U_{\ca m}(z,s,k) \in C^\infty(\G\backslash\H,k)$.
It should be clear from
the context whether we mean this new notion of weight or the
previous definition of weight. Trivially, if $\psi \in
C^\infty(\G\backslash\H,k)$ then $|\psi|$ has weight zero (in
either definition). We define the \it Maass raising and
lowering operators \rm by
$$
R_k=2iy\frac{d}{dz} +\frac{k}{2}, \ L_k=-2iy\frac{d}{d\bar z}
-\frac{k}{2}.
$$
It is an elementary exercise to show that
$$
R_k: C^\infty(\G\backslash\H,k)\rightarrow C^\infty(\G\backslash\H,k+2),\ 
L_k:
C^\infty(\G\backslash\H,k)\rightarrow C^\infty(\G\backslash\H,k-2).
$$
For $n>0,$ we write $R^n$ for $R_{k+2n-2} \cdots R_{k+2} R_k$ and 
$L^n$ for $L_{k-2n+2} \cdots L_{k-2}L_k$. (To simplify the notation we 
omit $k$ from the notation of the operators $L^n$ and $R^n$. It will 
usually be clear in each case). We also let $L^0$ and $R^0$ be the 
identity operator.

The hyperbolic Laplacian $\Delta$ can be realized as
$$
\Delta =-L_{2}R_0  = -R_{-2}L_0.  \tag 8.12
$$
By a direct calculation (see also Lemma 9.2 of \cite{JO})
$$
\align
R_k U_{\ca m}(z,s,k)&=(s+k/2)U_{\ca m}(z,s,k+2)-4\pi m U_{\ca
m}(z,s+1,k+2)\tag 8.13 \\
L_k U_{\ca m}(z,s,k)&=(s-k/2)U_{\ca m}(z,s,k-2). \tag 8.14
\endalign
$$
To see what happens to the Fourier expansions when we raise or lower the weight we need the next lemma. Set
$$
\omega (n,m,i,j) = (-4 \pi m)^i \left( \frac{m}{|m|} \right)^j 
\frac{(2n)!}{(i+j)! (i-j)! (n-i)!}.
$$

\proclaim{Lemma 8.1} For $n \geqslant 0$ we have
$$
\align
R^n\bigl(W_s(mz)\bigr) & = \sum_{i=0}^n \sum_{j=-i}^i \omega (n,m,i,j) y^i W_{s+j}(mz), \tag 8.15 \\
L^n\bigl(W_s(mz)\bigr) & = \sum_{i=0}^n \sum_{j=-i}^i \omega (n,-m,i,j) y^i W_{s+j}(mz), \tag 8.16 \\
R^n(y^s) = L^n(y^s) & = s(s+1) \cdots (s+n-1)y^s. \tag 8.17
\endalign
$$
\p See Sections 4 and 5 of \cite{O2} for the proofs of $(8.15)$ and $(8.16)$. A simple computation gives $(8.17)$. \bs

It follows from this lemma that
$$
R^n\bigl(W_s(mz)\bigr), L^n\bigl(W_s(mz)\bigr) \ll \sum_{i=0}^n \sum_{j=-i}^i (|m| y)^i |W_{s+j}(mz)|.
$$
Combine this with $(8.6)$ to prove that for any
$s$ with $1/2 \leqslant $ Re$(s) \leqslant 1-\delta$, $\delta >0$ and any $l \geqslant 
0$,
$$
R^n\bigl(W_{s}(mz)\bigr), L^n\bigl(W_{s}(mz)\bigr) \ll (|m|y)^{n-2l-3/2}(|s|^{2l+2+n} +1) |\G(s)| \tag 8.18
$$
as $y \rightarrow \ci$. The implied constant depends on $n$, $l$ and $\delta$.

\proclaim{Lemma 8.2} For $n \geqslant 0$ we have
$$
R^n\bigl(\eta_j(z)\bigr), L^n\bigl(\eta_j(z)\bigr) \ll (|t_j|^n +1) 
y_\G(z)^{1/2}+ (|t_j|^{2n+5} +1) y_\G(z)^{-3/2} 
$$
with the implied constant depending on $n$ and $\G$ alone.
\p
First consider 
$$
\rho_{\ca j}(0)R^n(y^{1-s_j})+\sum_{m\neq 0}
\rho_{\ca j}(m)R^n(W_{s_j}(mz)).\tag 8.19
$$
With $(8.4)$, $(8.5)$, $(8.7)$, $(8.17)$ and $(8.18)$ we see that 
$(8.19)$ is uniformly convergent in $z$ and bounded by a constant times
$$
(|t_j|^n +1) y^{1/2}+ (|t_j|^{2n+5} +1) y^{-3/2}
$$
as $y \rightarrow \ci$. Thus $(8.19)$ must equal $R^n(\eta_j(\sa z))$. 
To get from this to the statement of the lemma it is easiest to introduce 
the following operator.
For $\tau \in PSL_2(\R)$, let the operator $\theta_{\tau, k}: C^\infty 
(\G \backslash \H ,k) \rightarrow C^\infty (\tau^{-1}\G \tau \backslash 
\H ,k)$ have the action
$$
\theta_{\tau, k} \psi(z) = \frac{\psi(\tau z)}{\ee(\tau, z)^k}. \tag 8.20
$$
We will need the easily verified fact that $\theta$ commutes with the 
raising and lowering operators:
$$
\align
\theta_{\tau, k-2} L_k &= L_k \theta_{\tau, k}, \tag 8.21\\
\theta_{\tau, k+2} R_k &= R_k \theta_{\tau, k}. \tag 8.22
\endalign
$$
Now we can say that
$$
\left| \left. \left( R^n \eta_j(w)\right)\right|_{w=\sa z} \right| = 
\left| \theta_{\sa, 2n} R^n \eta_j(z) \right|
= \left| R^n \theta_{\sa, 0}  \eta_j(z) \right| = \left| 
R^n(\eta_j(\sa z)) \right|.
$$
Similarly for the lowering operator $L^n$ and the lemma follows. \bs

\proclaim{Lemma 8.3} For $n \geqslant 0$ we have
$$
\int_T^{T+1}  \left| R^n E_{\ca}( z, 1/2+ir)\right|^2 \, dr, 
\ \ \int_T^{T+1}  \left| L^n E_{\ca}( z, 1/2+ir)\right|^2 \, dr \ll 
T^{4n+12} y_\G(z).
$$
\p
With (4.6) and(8.18) we see that, for $r$ in $[T,T+1]$ and any integer $l \geqslant 0$,
$$
R^n E_{\ca}(\sb z, 1/2+ir) \ll T^{n} y^{1/2}+ \sum_{m \neq 0} \left| \phi_{\ca \cb}(m,1/2+ir)\right| 
(|m|y)^{n-2l-3/2} T^{2l+2+n} e^{-\pi T/2},
$$
so that
$$
\multline
\int_T^{T+1}  \left| R^n E_{\ca}(\sb z, 1/2+ir)\right|^2 \, dr \ll T^{2n} y 
\\
+  T^{2n+2l+2} e^{-\pi T/2} y^{n-2l-1}
\int_T^{T+1} \sum_{m \neq 0} \left| \phi_{\ca \cb}(m,1/2+ir)\right| 
|m|^{n-2l-3/2}  \, dr\\
+  T^{2n+4l+4} e^{-\pi T} y^{2n-4l-3}
\int_T^{T+1} \sum_{m_1, m_2 \neq 0} \left| \phi_{\ca \cb}(m_1,1/2+ir)\phi_{\ca \cb}(m_2,1/2+ir)\right| 
|m_1 m_2|^{n-2l-3/2}  \, dr .
\endmultline
\tag 8.23
$$
With (8.11) it is clear that the sums in the integrals on the right side 
of (8.23) 
are absolutely convergent for $l=\lceil n/2 +1 \rceil$ and uniformly 
bounded for $T \leqslant r \leqslant T+1$. Interchanging sums and 
integrals in (8.23) is now justified with Corollary 8.6 below and 
with (8.10) and the Cauchy-Schwartz inequality we obtain the lemma for the raising operator. The argument for $L^n$ is identical. \bs

We shall have frequent need of the next three standard analysis results.

\proclaim{Theorem 8.4} Let $z$ be in a set $S \subseteq \C$. Suppose the functions $f_k(z)$ are smooth and $\sum_{k=1}^\infty f_k(z)$ converges pointwise on $S$. If $\sum_{k=1}^\infty \frac{d}{dz}f_k(z)$ converges uniformly then
$$
\frac{d}{dz} \sum_{k=1}^\infty f_k(z)=\sum_{k=1}^\infty \frac{d}{dz}f_k(z).
$$
\endproclaim

\proclaim{Theorem 8.5} Suppose the functions $g_k(r)$ are integrable on $[a,b]$ and $\lim_{k \to \infty} g_k(r)$ converges pointwise then
$$
\lim_{k \to \infty} \int_a^b g_k(r) \, dr=  \int_a^b \lim_{k \to \infty} g_k(r) \, dr
$$
if, for all $k$, $|g_k(r)| \leqslant C$ for some fixed constant $C$.
\endproclaim

\proclaim{Corollary 8.6} Suppose the functions $h_k(r)$ are smooth on $[a,b]$ and $\sum_{k=1}^\infty |h_k(r)|$ exists and is uniformly bounded on $[a,b]$ then
$$
\sum_{k=1}^\infty \int_a^b h_k(r) \, dr=  \int_a^b \left( \sum_{k=1}^\infty h_k(r)\right) \, dr.
$$
\endproclaim

Theorem 8.4 is a weak form of Theorem 7.17 in \cite{Ru}. Theorem 8.5 follows from Lebesgue's Dominated Convergence
Theorem as in  \cite{Ru}, Theorem 10.32. The corollary follows directly from Theorem 8.5.

\newpage

$$\text{\bf 9. Proof of Proposition A}$$

\proclaim{Proposition A} For $k \in 2\Z$ and $\sigma=\Re (s) >1$ the series
$U_{\ca m}(z,s,k)$, $Q_{\ca m}(z,s,1;\overline{f})$, $Q'_{\ca m}(z,s,1;\overline{f})$ and $G_{\ca m}(z,s-1;\overline{F})$ converge absolutely and uniformly on compact sets to analytic functions of $s$. For these $s$ we
have
$$
\align
U_{\ca 0}(z,s,k) & \ll   y_\G(z)^{\sigma},\tag i\\
U_{\ca m}(z,s,k) & \ll  1, \ \ \ m>0,\tag ii\\
Q_{\ca m}(z,s,1;\overline{f}) & \ll   y_\G(z)^{1/2-\sigma/2}, \ \ \ m \geqslant 0,\tag iii\\
yQ'_{\ca m}(z,s,1;\overline{f}) & \ll   (|m|+1) y_\G(z)^{1/2-\sigma/2}, \ \ \ m \geqslant 0,\tag iv\\
yG_{\ca m}(z,s-1;\overline{F}) & \ll  y_\G(z)^{1/2-\sigma/2}, \ \ \ m \geqslant 0 \tag v
\endalign
$$
where the implied constants depends on $s$, $k$, $f$ and $\G$ but not $m$.
\p
The Eisenstein series $E_\ca(z,s)$ given by (4.5) equals $U_{\ca 0}(z,s,0)$, is absolutely convergent for $\Re (s)>1$ and satisfies
$$
E_{\ca}(\sb z,s)=\delta_{\ca \cb}y^s +\phi_{\ca \cb}(s)y^{1-s}+O(e^{-2\pi y}) \tag 9.1
$$
as $y \to \infty$ by (4.7). Hence
$$
U_{\ca 0}(z,s,k) \ll E_\ca(z,\sigma) \ll y_\G(z)^{\sigma},
$$
which is (i).

\vskip .10in
For $m>0$
we have (with $|e(m \sa^{-1} \g z)| \leqslant 1$)
$$
\align U_{\ca m}(\sa z,s,k) &\ll y^\sigma e^{-2\pi
my}+\sum \Sb \g \in \G_\ca\backslash\G
\\ \g \neq identity \endSb \Im(\sa^{-1}\g\sa z)^\sigma \\
&\ll y^\sigma e^{-2\pi my}+|E_\ca(\sa z,\sigma)-y^\sigma| \ll 1.
\endalign
$$
At any other cusp $\cb \neq \ca$
$$
U_{\ca m}(\sb z,s,k)\ll E_\ca(\sb z,\sigma)\ll \phi_{\ca \cb}(s)y^{1-\sigma} \ll 1.
$$
We have shown statement (ii).

\vskip .10in
By Lemma 4.1 and subsequent discussion we have, for any $\epsilon>0$, 
$$
F_\ca(\sb z) \ll y^\epsilon+y^{-\epsilon}+1 \tag 9.2
$$
for all $z$ in $\H$. It  is then apparent that
$$
F_\ca(\g \sb z)=F_\ca(\sa \sa^{-1} \g \sb z) \ll \Im (\sa^{-1} \g \sb z)^\epsilon + \Im (\sa^{-1} \g \sb z)^{-\epsilon}+1
$$
for any cusp $\cb$ and any $z$ in $\H$. The implied constant depends solely on $\epsilon, f$ and $\G$. In the case $\ca =\cb$ we may improve $(9.2)$. The Fourier expansion of $F_\ca$ yields
$$
F_\ca(\sb z) = \int_\ca^\cb f(w)\, dw +\frac{1}{2\pi i} \sum_{n=1}^\infty \frac{a_\cb(n)}{n} e(nz) \tag 9.3
$$
with $a_\cb(n)$ the $n$th Fourier coefficient of $f$ at the cusp $\cb$. Thus, when $\ca = \cb$,
$$
F_\ca(\sa z) \ll e^{-2\pi y} \text{ \ as } y \rightarrow \infty. \tag 9.4
$$
Consequently
$$
\align
Q_{\ca m}(\sa z, s,1; \overline{f}) & \ll \sum_{\g \in \G_\ca \backslash \G}  |F_{\ca}(\g \sa  z)| \Im(\sa^{-1} \g \sa z)^\sigma \\
& \ll e^{-2\pi y} y^\sigma + \sum \Sb \g \in \G_\ca \backslash \G \\ \g \neq identity \endSb  \left( \Im(\sa^{-1} \g \sa z)^{\sigma+\epsilon} + \Im(\sa^{-1} \g \sa z)^{\sigma-\epsilon} + \Im(\sa^{-1} \g \sa z)^\sigma \right)\\
&\ll y^{1-\sigma +\epsilon}
\endalign
$$
for $\sigma> 1+\epsilon $ as $y \rightarrow \infty$  by (9.1).
When $\ca \neq \cb$ we do not need to worry about the $y^s$ term and
$$
\align
Q_{\ca m}(\sb z, s,1; \overline{f}) & \ll E_\ca(\sb z, \sigma+\epsilon) + E_\ca(\sb z, \sigma-\epsilon)+ E_\ca(\sb z, \sigma)\\
&\ll y^{1-\sigma +\epsilon}
\endalign
$$
for $\sigma >1+\epsilon$ as $y \rightarrow \infty$. 
Choose $\epsilon =(\sigma -1)/2$ for simplicity and we have demonstrated that 
$$
Q_{\ca m}(z, s,1; \overline{f}) \ll y_\G(z)^{1/2-\sigma/2}
$$
for $\sigma > 1$ and an implied constant depending on $\sigma, f$ and 
$\G$ alone. This is (iii).

\vskip .10in
Taking derivatives we see, for any $\g$ in $PSL_2(\R)$,
$$
\multline
2iy \frac{d}{dz} \left(\overline{F_{\ca}(\sa \g z)} \Im( \g z)^s e(m  \g 
z)\right)
=s\overline{F_{\ca}(\sa \g z)} \Im( \g z)^s e(m  \g z)\ee(\g, z)^{-2} \\
-4\pi m \overline{F_{\ca}(\sa \g z)} \Im( \g z)^{s+1} e(m  \g z)\ee(\g, 
z)^{-2}.
\endmultline
$$
Hence
$$
\align
yQ'_{\ca m}(z, s,1; \overline{f}) & \ll |s|\sum_{\g \in \G_\ca \backslash 
\G}  |F_{\ca}(\g   z)| \Im(\sa^{-1} \g  z)^\sigma + |m|\sum_{\g \in 
\G_\ca \backslash \G}  |F_{\ca}(\g   z)| \Im(\sa^{-1} \g  z)^{\sigma+1}
\\
&\ll (|m|+1)y_\G(z)^{1-\sigma +\epsilon}
\endalign
$$
and set $\epsilon =(\sigma -1)/2$ as before to obtain (iv).

\vskip .10in
Finally, it is easy to see that
$$
yG_{\ca m}(z,s;\overline{F})=\sum_{\g \in \G_\ca \backslash \G} 
\overline{F_\ca(\g z)}  \Im(\sa^{-1} \g z)^{s+1} e(m 
\sa^{-1} \g z)\ee(\sa^{-1} \g, z)^{-2}.
$$
So the argument used for $Q_{\ca m}(z, s,1; \overline{f})$ applies to $G_{\ca m}(z,s-1;\overline{F})$ yielding (v) and the proof is complete. \bs

$$\text{\bf 10. Proof of Propositions B and C}$$

\proclaim{Proposition B} For $k \in 2\Z$  the (Eisenstein) series
$U_{\ca 0}(z,s,k)$  has a meromorphic
continuation to all $s$ with $\Re (s) >1-\delta_\G$. We
have
$$
U_{\ca 0}(z,s,k) \ll   y_\G(z)^{\sigma}.
$$
for these $s$ with the implied constant depending on $s$, $k$ and $\G$.
The only possible pole in this region appears at $s=1$
when $k=0$. It is a simple pole with residue $1/V$.
\p
The Eisenstein series $E_\ca(z,s)=U_{\ca 0}(z,s,0)$ has a meromorphic continuation to all $s$ in $\C$. This is shown in Chapter 6 of \cite{I1}. There it is also shown that the Fourier expansion 
$$
E_{\ca}(\sb z,s)=\delta_{\ca \cb}y^s +\phi_{\ca \cb}(s)y^{1-s}+\sum_{m \neq 0} \phi_{\ca \cb}(m,s)W_s(mz) \tag 10.1
$$
is valid for all $s$ in $\C$ except at the poles of $E_\ca(z,s)$. In particular, for $\Re (s) >1-\delta_\G$ the expansion (10.1) is valid except at $s=1$ where $\phi_{\ca \cb}(s)$ has a simple pole. Proposition 6.13 of \cite{I1} shows that the residue of $E_\ca(z,s)$ at $s=1$ is $1/V$. The coefficients $\phi_{\ca \cb}(m,s)$ are analytic for $\sigma=\Re (s) >1-\delta_\G$ and, for all $s$, satisfy
$$
\phi_{\ca \cb}(m,s) \ll |m|^\sigma+|m|^{1-\sigma} \tag 10.2
$$
with an implied constant depending uniformly on $s$ (away from poles) and $\G$. This result is stated in \cite{I1} (6.19) and proved in \cite{JO} Proposition 7.2.

\vskip .10in
Now by (8.13), (8.14)
$$
\align
R^n E_\ca(z,s) & =s(s+1) \cdots (s+n-1) U_{\ca 0}(z,s,2n),\\
L^n E_\ca(z,s) & =s(s+1) \cdots (s+n-1) U_{\ca 0}(z,s,-2n),
\endalign
$$
so we may obtain the meromorphic continuation of $U_{\ca 0}(z,s,k)$ from $E_\ca(z,s)$.
For example, with Lemma 8.1, Theorem 8.4, (8.18), (10.1) and (10.2) we have
$$
\align
R^n E_\ca(\sb z,s) & =\delta_{\ca \cb}R^n y^s +\phi_{\ca \cb}(s)R^n y^{1-s}+R^n \sum_{m \neq 0} \phi_{\ca \cb}(m,s)W_s(mz),\\
 & =\delta_{\ca \cb}s(s+1) \cdots (s+n-1) y^s +\phi_{\ca \cb}(s)(1-s)(2-s) \cdots (n-s) y^{1-s}\\
& \qquad + \sum_{m \neq 0} \phi_{\ca \cb}(m,s)R^n W_s(mz).
\endalign
$$
For $n>0$ the pole of $\phi_{\ca \cb}(s)$ at 1 is eliminated and we see that 
$$
R^n E_\ca(\sb z,s) \ll \delta_{\ca \cb}y^\sigma + y^{1-\sigma}.
$$
Similarly for $L^n$ and we have shown that
$$
U_{\ca 0}(\sb z,s,k) \ll \delta_{\ca \cb}y^\sigma + y^{1-\sigma} \tag 10.3
$$
as $y \to \infty$ for an implied constant depending on $s$, $k$ and $\G$. \bs

\proclaim{Proposition C} For $k \in 2\Z$ and $m > 0$ the Poincar\'e series
$U_{\ca m}(z,s,k)$  has a
continuation to an analytic function for all $s$ with $\Re (s) >1-\delta_\G$. We
have
$$
U_{\ca m}(z,s,k) \ll   y_\G(z)^{1/2}
$$
for these $s$ with the implied constant depending on $s$, $m$, $k$ and $\G$.
\p
We first look at the case $k=0$. 
In the spectral decomposition (8.1) of $U_{\ca m}(z,s,0)$ the inner 
products may be found explicitly as in 
 $(8.4)$ of \cite{JO} and \cite{I3}, Chapter 17:
$$
\multline
U_{\ca m}(z,s,0)\pi^{-1/2}(4\pi m)^{s-1/2}\G(s)=\sum_{j=1}^\infty \G(s-s_j)
\G(s-1+s_j)\overline{\rho_{\ca j}}(m)\eta_j(z)\\
+\frac{1}{4\pi}\sum_\cb \int_{-\infty}^\infty \G(s-1/2-ir)\G(s-1/2+ir)
\overline{\phi_{\ca \cb}}(m,1/2+ir)E_\cb(z,1/2+ir)\, dr
\endmultline \tag 10.4
$$
where, referring to $(8.3)$, the $\rho_{\ca j}(m)$ are the Fourier coefficients of the Maass forms $\eta_j$. Let
$$
U_{\ca m}(z,s,0)_{DISC}=\sum_{j=1}^\infty \G(s-s_j)
\G(s-1+s_j)\overline{\rho_{\ca j}}(m)\eta_j(z),
$$
the discrete spectral component. We shall examine this first. 
Use $(8.5)$ and $(8.7)$ to get
$$
\G(s-s_j)\G(s-1+s_j)\overline{\rho_{\ca j}}(m) \ll \frac{|t_j|^{2\sigma - 1/2}}{\sqrt{|m|}}e^{-\pi |t_j|/2}. \tag 10.5
$$

\proclaim{Lemma 10.1} Let $m$ and $n$ be integers with $m>0$, $n 
\geqslant 0$. As functions of $s$, $R^n\bigl(U_{\ca 
m}(z,s,0)_{DISC}\bigr)$ and $L^n\bigl(U_{\ca m}(z,s,0)_{DISC}\bigr)$ are 
analytic for $\Re (s) > 1-\delta_\G$. For these $s$ and an implied 
constant depending on $s$, $n$ and $\G$ alone we have
$$
R^n\bigl(U_{\ca m}(z,s,0)_{DISC} \bigr), \, \, L^n\bigl(U_{\ca 
m}(z,s,0)_{DISC} \bigr)  \ll |m|^{-1/2} y_\G (z)^{1/2}.
$$
\p Let 
$$
J(z,s)=\sum_{j=1}^\infty \G(s-s_j)
\G(s-1+s_j)\overline{\rho_{\ca j}}(m)R^n\bigl(\eta_j(z) \bigr).
$$
 With $(10.5)$, Lemma 8.2 and $(8.2)$ it follows that, for fixed $s$, the series $J(z,s)$ converges uniformly for  $z$ in any compact set, say, and is bounded by $|m|^{-1/2} y_\G (z)^{1/2}$. Hence, with Theorem 8.4,
$$
R^n\bigl(U_{\ca m}(z,s,0)_{DISC} \bigr) = J(z,s).
$$ 
We also see that $J(z,s)$ converges uniformly for $s$ in compact sets with $\Re (s) > 1-\delta_\G$ giving an analytic function of $s$.
Similarly for $L^n$. \bs

To deal with the continuous spectral component,
$$
U_{\ca m}(z,s,0)_{CONT}=\frac{1}{4\pi}\sum_\cb \int_{-\infty}^\infty \G(s-1/2-ir)\G(s-1/2+ir)
\overline{\phi_{\ca \cb}}(m,1/2+ir)E_\cb(z,1/2+ir)\, dr,
$$
we shall need the next lemma.

\proclaim{Lemma 10.2} For $\psi(r)$ smooth on $[T,T+1]$ we have
$$
\frac{d}{dz} \int_T^{T+1} \psi(r) E_{\ca}(\sb z, 1/2+ir) \, dr
=  \int_T^{T+1} \psi(r) \left( \frac{d}{dz} E_{\ca}(\sb z, 1/2+ir)\right) \, dr.
$$
\p
With the Fourier expansion (4.6) we have
$$
E_{\ca}(\sb z,1/2+ir)=\delta_{\ca \cb}y^{1/2+ir} +\phi_{\ca \cb}(1/2+ir)y^{1/2-ir}+\sum_{m \neq 0} \phi_{\ca \cb}(m,1/2+ir)W_{1/2+ir}(mz).
$$
 Combine (8.6), (8.9) and (8.11) to see that
$$
\int_T^{T+1} \left( \delta_{\ca \cb} \left|y^{1/2+ir}\right| +\left|\phi_{\ca \cb}(1/2+ir)y^{1/2-ir}\right|+\sum_{m \neq 0} \left|\phi_{\ca \cb}(m,1/2+ir)W_{1/2+ir}(mz)\right|\right) \,dr < \infty,
$$
and hence, by  Corollary 8.6,
$$
\multline
\int_T^{T+1} \psi(r) E_{\ca}(\sb z, 1/2+ir) \, dr
=\int_T^{T+1} \left( \delta_{\ca \cb} y^{1/2+ir} + \phi_{\ca \cb}(1/2+ir)y^{1/2-ir}\right) \, dr
\\
+\sum_{m \neq 0} \int_T^{T+1} \phi_{\ca \cb}(m,1/2+ir)W_{1/2+ir}(mz) \,dr 
\endmultline
\tag 10.6
$$
and similarly, using (8.15),
$$
\multline
\int_T^{T+1} \psi(r) \left( \frac{d}{dz} E_{\ca}(\sb z, 1/2+ir)\right) \, dr
=\int_T^{T+1} \left( \delta_{\ca \cb} \frac{d (y^{1/2+ir})}{dz} 
+  \phi_{\ca \cb}(1/2+ir) \frac{d (y^{1/2-ir})}{dz} \right)\, dr \\
+\sum_{m \neq 0} \int_T^{T+1} \phi_{\ca \cb}(1/2+ir,m)\left( \frac{d}{dz} W_{1/2+ir}(mz)\right) \,dr 
\endmultline
\tag 10.7
$$
To demonstrate that the derivative of (10.6) equals (10.7) we need to show that all the corresponding components are equal. For example
$$
\frac{d}{dz} \int_T^{T+1} \psi(r)  y^{\sigma+ir} \, dr
=  \int_T^{T+1} \psi(r)  \frac{d (y^{\sigma+ir})}{dz} \, dr \tag 10.8
$$
because, using  Theorem 8.5, it is easy to check that, for each $y$, 
$$
((y+h)^{\sigma+ir}-y^{\sigma+ir})/h
$$ is uniformly bounded for $r$ in $[T,T+1]$ as $h \to 0$.  
In the same way
$$
\frac{d}{dz} \int_T^{T+1} \psi(r)  W_{1/2+ir}(mz) \, dr
=  \int_T^{T+1} \psi(r)  \left( \frac{d}{dz}W_{1/2+ir}(mz)\right) \, dr,
$$
completing the proof. \bs

Moreover, the same arguments and Lemma 8.1 imply
$$
\align
R^n \int_T^{T+1} \psi(r) E_{\ca}(\sb z, 1/2+ir) \, dr
=  \int_T^{T+1} \psi(r) \left( R^n E_{\ca}(\sb z, 1/2+ir)\right) \, dr, 
\tag 10.9\\
L^n \int_T^{T+1} \psi(r) E_{\ca}(\sb z, 1/2+ir) \, dr
=  \int_T^{T+1} \psi(r) \left( L^n E_{\ca}(\sb z, 1/2+ir)\right) \, dr. \tag 10.10
\endalign
$$

Returning to our continuous spectral component,
$$
R^n U_{\ca m}(z,s,0)_{CONT}=\frac{1}{4\pi}\sum_\cb R^n \int_{-\infty}^\infty \G(s-1/2-ir)\G(s-1/2+ir)
\overline{\phi_{\ca \cb}}(m,1/2+ir)E_\cb(z,1/2+ir)\, dr.
$$
If we restrict our attention to $r$ in $[T,T+1]$ we find
$$
\align
\frac{1}{4\pi}\sum_\cb R^n & \int_{T}^{T+1} \G(s-1/2-ir)\G(s-1/2+ir)
\overline{\phi_{\ca \cb}}(m,1/2+ir)E_\cb(z,1/2+ir)\, dr \\
& = \frac{1}{4\pi}\sum_\cb  \int_{T}^{T+1} \G(s-1/2-ir)\G(s-1/2+ir)
\overline{\phi_{\ca \cb}}(m,1/2+ir)R^n E_\cb(z,1/2+ir)\, dr\\
& \ll 
T e^{-\pi T}\sqrt{ \int_{T}^{T+1} \left| \phi_{\ca \cb}(m,1/2+ir)\right|^2 \, dr}
\sqrt{ \int_{T}^{T+1} \left| R^n E_\cb(z,1/2+ir) \right|^2 \, dr}\\
& \ll 
T^{2n +8} e^{-\pi T/2}|m|^{-1/2} y_\G(z)^{1/2},
\endalign
$$
where we used (10.9) to get the second line, the Cauchy-Schwartz inequality and (8.7) for line three, and (8.10) and Lemma 8.3 for the last line.
Therefore, repeating the argument for $L^n$, we have shown the following.

\proclaim{Lemma 10.3} Let $m$ and $n$ be integers with $m>0$, $n 
\geqslant 0$. As functions of $s$, $R^n\bigl(U_{\ca 
m}(z,s,0)_{CONT}\bigr)$ and $L^n\bigl(U_{\ca m}(z,s,0)_{CONT}\bigr)$ are 
analytic for $\Re (s) > 1-\delta_\G$. For these $s$ and an implied 
constant depending on $s$, $n$ and $\G$ alone we have
$$ R^n\bigl(U_{\ca m}(z,s,0)_{CONT} \bigr), \, \, L^n\bigl(U_{\ca 
m}(z,s,0)_{CONT} \bigr)  \ll |m|^{-1/2} y_\G (z)^{1/2}.
$$
\endproclaim

\vskip .10in
We may now finish the proof of Proposition C.
With (8.13) we see that
$$
\align
R_0 U_{\ca m}(z,s,0)& =sU_{\ca m}(z,s,2)-4\pi m U_{\ca
m}(z,s+1,2),\\
R_2 R_0 U_{\ca m}(z,s,0)& =s(s+1)U_{\ca m}(z,s,4)-4\pi m (2s+2)U_{\ca
m}(z,s+1,4)+(4\pi m)^2 U_{\ca
m}(z,s+2,4).
\endalign
$$
In general, for $k \geqslant 0$,
$$
\multline
U_{\ca m}(z,s,2k)=\frac{1}{s(s+1) \cdots (s+k-1)} \left( R^k U_{\ca m}(z,s,0)\right. \\
\left. + p_1(m,s)U_{\ca m}(z,s+1,2k)+ \cdots +p_k(m,s)U_{\ca m}(z,s+k,2k)\right) 
\endmultline 
\tag 10.11
$$
with polynomials $p_i$ in $m$ and $s$. Therefore, using Lemmas 10.1, 10.3 and Proposition A part (ii), the right side of (10.11) is analytic for $\Re(s)>1-\delta_\G$ and bounded by $y_\G (z)^{1/2}$.
Similarly for $k<0$. \bs

$$\text{\bf 11. Proof of Proposition D}$$

\proclaim{Proposition D} For $m \geqslant 0$, both series $(s-1)Q_{\ca 
m}(z,s,1;\overline{f})$ and $Q'_{\ca m}(z,s,1;\overline{f})$ have  
continuations to analytic functions of $s$ with 
$\Re (s) >1-\delta_\G$. For these $s$ values
$$
(s-1)Q_{\ca m}(z,s,1;\overline{f}), \ \ \ yQ'_{\ca m}(z,s,1;\overline{f}) \ll  
y_\G(z)^{1/2}.
$$
The implied constant depends on $s$, $m$, $f$ and $\G$. Also $Q_{\ca 
m}(z,s,1;\overline{f})$ has a simple pole at $s=1$ with residue 
$2i\overline{\s{f}{P_{\ca m}(\cdot)_2}}$.
\p
By Proposition A part (iii), $Q_{\ca m}(z,s,1;\overline{f})$ is certainly square integrable for $\Re(s)>1$. In other words 
$$
\langle Q_{\ca m}(\cdot,s,1;\overline{f}),Q_{\ca m}(\cdot,s,1;\overline{f})\rangle < \infty.
$$ 
The
Roelcke-Selberg decomposition, (8.1), yields
$$
\multline
Q_{\ca m}(z,s,1;\overline{f})=\sum_{j=0}^\infty\langle Q_{\ca m}(\cdot,s,1;\overline{f}),\eta_j\rangle
\eta_j(z)\\
+\frac{1}{4\pi }\sum_\cb \int_{-\infty}^{\infty}\langle
Q_{\ca m}(\cdot,s,1;\overline{f}) ,E_\cb(\cdot,1/2+ir)\rangle E_\cb(z,1/2+ir)\,dr.
\endmultline \tag 11.1
$$
To understand the inner products appearing in (11.1) we make use of the next lemma.

\proclaim{Lemma 11.1} Let $\xi_1$, $\xi_2$ and $\psi$ be any smooth $\G$ invariant functions (not necessarily in $L^2(\G \backslash \H)$). If $(\Delta - \lambda) \xi_1 =\xi_2$, $(\Delta - \lambda') \psi =0$ and
$$
\align
\xi_1, R_0\xi_1, \Delta \xi_1 & \ll y_\G(z)^A,
\\ \psi, R_0 \psi & \ll y_\G(z)^B
\endalign
$$
for $A+B<0$ and $R_0=2iy \frac{d}{dz}$ the raising operator, then
$$
\langle \xi_1,\psi \rangle = \frac 1{\lambda' -\lambda} \langle \xi_2,\psi \rangle.
$$
\p
We simply have
$$
\langle \xi_1,\psi \rangle = \frac 1{\lambda' -\lambda} \langle \xi_1,(\Delta - \overline{\lambda})\psi \rangle 
= \frac 1{\lambda' -\lambda} \langle (\Delta - \lambda)\xi_1,\psi \rangle 
= \frac 1{\lambda' -\lambda} \langle \xi_2,\psi \rangle.
$$
To justify switching $\Delta$ from the right side of the inner product to the left side requires the growth assumptions we stated. See \cite{JO} Proposition 9.3 and Corollary 9.4 for the proof of this. \bs

Now for all $n \in \Z$,
$$
\align
(\Delta - s(1-s))Q_{\ca m}(z,s,n;\overline{f})  = & -8 \pi i m Q_{\ca m}(z,s+2,n-1;\overline{f})\\
&+ 4\pi ms Q_{\ca m}(z,s+1,n;\overline{f})\\
&+ 2i s Q_{\ca m}(z,s+1,n-1;\overline{f}).
\endalign
$$
We want to apply this lemma to $\xi_1 =Q_{\ca m}(z,s,n;\overline{f})$ and $\xi_1 = \eta_j$, (recall that $(\Delta -s_j(1-s_j))\eta_j =0$). To check the growth conditions we will need the following result.

\proclaim{Proposition E} For $-n \leqslant 0$ the  series $Q_{\ca m}(z,s+n+1,-n;\overline{f})$ is an analytic function of $s$ for $\Re (s) >1-\delta_\G$. Also for these $s$  we have 
$$
Q_{\ca m}(z,s+n+1,-n;\overline{f}), \ \ R_0 Q_{\ca m}(z,s+n+1,-n;\overline{f}) \ll  e^{-\pi y_\G(z)} 
$$
with the implied constant depending on $n,m,f,s$ and $\G$ alone. 
\endproclaim

The proof of this proposition follows at the end of this section and depends on the nice fact that
$Q_{\ca m}(z,s+n+1,-n;\overline{f})$ with $-n \leqslant 0$ can be expressed as a linear combination of Poincar\'e series $U_{\ca m}(z,s,k)$ multiplied by something with exponential decay at the cusps (that is, $\overline{f}$ with its weight lowered by $L$).

\vskip .10in
Now we have $\eta_j(z)$, $R_0\eta_j(z) \ll y_\G(z)^{1/2}$ by Lemma 8.2 and 
$$
Q_{\ca m}(z,s,1;\overline{f}), R_0Q_{\ca m}(z,s,1;\overline{f}), \Delta Q_{\ca m}(z,s,1;\overline{f}) \ll y_\G(z)^{1/2-\sigma/2}
$$
for $\sigma= \Re(s)>1$ by Proposition A, parts (iii) and (iv) and Proposition E. So we may use Lemma 11.1 to get, for $\Re(s)>2$,
$$
\multline
\langle Q_{\ca m}(\cdot,s,1;\overline{f}),\eta_j\rangle = \frac{1}{(s_j-s)(1-s_j-s)}\Big( 
-8 \pi i m \langle Q_{\ca m}(\cdot,s+2,0;\overline{f}),\eta_j\rangle\\
+ 4\pi ms \langle Q_{\ca m}(\cdot,s+1,1;\overline{f}),\eta_j\rangle
+ 2i s \langle Q_{\ca m}(\cdot,s+1,0;\overline{f}),\eta_j\rangle \Big).
\endmultline
$$
We can repeat this procedure $W$ times in all to obtain, again for $\Re(s)>2$,
$$
\langle Q_{\ca m}(\cdot,s,1;\overline{f}),\eta_j\rangle = \sum_l \frac{P_l(m,s)}{R_l(s_j,s)} 
\langle Q_{\ca m}(\cdot,s+W+c_l,1-d_l;\overline{f}),\eta_j\rangle, \tag 11.2
$$
with integers $c_l,d_l$ satisfying $0 \leqslant c_l,d_l \leqslant W$, 
$d_l \le W+c_l$, $P_l(m,s)$ a polynomial in $m,s$ alone of degree $W$ in 
$m$ and of 
degree $W$ in $s$ and $R_l(s_j,s)$ a polynomial in $s_j,s$ alone of degree $2W$ in $s_j$ and of degree $2W$ in $s$. In fact
$$
R_l(s_j,s)=\prod_b (s_j-b-s)(1-s_j-b-s) \tag 11.3
$$
where, for each $l$, the product is over some subset of integers $b$ in $\{0,1, \cdots ,2W\}$ of cardinality $W$.

\vskip .10in
The finite sum on the right of (11.2) may be used to give the analytic continuation of the inner product on the left and to bound it.
For our purposes we are only interested in getting the analytic continuation to $\Re (s)>1-\delta_\G$. Examining each term on the right of $(11.2)$ we see that if $d_l=0$ then we have
$$
Q_{\ca m}(z, s+W+c_l,1;\overline{f}) \ll y_\G(z)^{1/4-W/2}\tag 11.4
$$
by Proposition A, (iii) for $W \geqslant 1$. Hence
$$
\langle Q_{\ca m}(\cdot,s+W+c_l,1;\overline{f}),\eta_j\rangle \ll \sqrt{ || y_\G(z)^{-1/4}|| \cdot || \eta_j ||} = \sqrt{ || y_\G(z)^{-1/4}|| } \ll 1. \tag 11.5
$$
For $0< d_l \leqslant W$,
Proposition E implies that
$$
Q_{\ca m}(z,s+W+c_l,1-d_l;\overline{f}) \ll  e^{-\pi y_\G(z)}. \tag 11.6
$$
Hence, as in $(8.6)$,
$$
\langle Q_{\ca m}(z,s+W+c_l,1-d_l;\overline{f}),\eta_j\rangle \ll  1. \tag 11.7
$$
Now combine  $(11.2)$, $(11.3)$, $(11.5)$ and $(11.7)$ to see that, for $j>0$,
$
\langle Q_{\ca m}(\cdot,s,1;\overline{f}),\eta_j\rangle
$
is an analytic function of $s$ for $\Re(s)>1 - \delta_\G$ and satisfies
$$
\langle Q_{\ca m}(\cdot,s,1;\overline{f}),\eta_j\rangle \ll |s_j|^{-2W} \ll |\lambda_j|^{-W} \tag 11.8
$$
for implied constants depending on $s,m,W,f$ and $\G$ alone.

\vskip .10in
For $j>0$ we can use (8.2), (8.8) and (11.8) to get
$$
\sum_{T\leqslant |\lambda_j|< T+1} \langle Q_{\ca m} 
(\cdot,s,1;\overline{f}),\eta_j\rangle \eta_j(z)
\ll T\left[y_\G(z)^{1/2} + T^{7/4} y_\G(z)^{-3/2}\right]T^{-W}.
$$
Therefore, (using any $W \geqslant 4$),
$$
\sum_{j=1}^\infty \langle Q_{\ca m}(\cdot,s,1;\overline{f}),\eta_j\rangle \eta_j(z)
\ll y_\G(z)^{1/2}  \tag 11.9
$$
for all $s$ with $\Re(s)>1-\delta_\G$ and an implied constant depending solely on $s$, $m$, $f$ and $\G$.

\vskip .10in
For $j=0$ the constant eigenfunction is $\eta_0=V^{-1/2}$. By unfolding we obtain
$$
\align
\langle Q_{\ca m}(\cdot,s,1;\overline{f}),\eta_0\rangle \eta_0 &=  \frac{-\overline{a_{\ca}(m)} \ \Gamma(s-1)}{2\pi i m (4\pi m)^{s-1}} \\
&= \frac{-\overline{a_{\ca}(m)} }{2\pi i m }\left( \frac{1}{s-1} +O(1)\right)\\
&= 2i \overline{\langle f, P_{\ca m}(\cdot )_2\rangle } \left( \frac{1}{s-1} +O(1)\right) \tag 11.10
\endalign
$$
as $s\rightarrow 1$ since $\langle f, P_{\ca m}(\cdot )_2\rangle = a_{\ca}(m)/(4\pi m)$ for
$$
f_\ca(z)=j(\sa, z)^{-2}f(\sa z)=\sum_{m=1}^\infty a_{\ca}(m) e(mz).
$$

\vskip .10in
With arguments similar  to those used for the discrete spectrum we now consider the continuous spectrum. For $P_l, R_l, c_l$ and $d_l$ identical to $(11.2)$,
$$
\langle Q_{\ca m}(\cdot,s,1;\overline{f}),E_\cb(\cdot, 1/2+ir)\rangle = \sum_l \frac{P_l(m,s)}{R_l(1/2+ir,s)} 
\langle Q_{\ca m}(\cdot,s+W+c_l,1-d_l;\overline{f}),E_\cb(\cdot, 1/2+ir)\rangle, \tag 11.11
$$
which is true for $\Re(s)> 2$ initially. 
 Here we employed (4.6), (8.6), (8.9) and (8.11) to get the bound
$$
E_{\ca}(z,1/2+ir) \operatornamewithlimits{\ll}  y_\G(z)^{1/2} \tag 11.12
$$
for $r \in [T,T+1]$ and an implied constant depending on $T$ and $\G$.

\vskip .10in
With $(11.4)$, $(11.6)$ and $(11.12)$ we see that (for $W\geqslant 1$) the right side of $(11.11)$ converges and gives the analytic continuation of the left side to $\Re(s)>1-\delta_\G$.
 Now 
$$
\multline
\int_T^{T+1} \langle Q_{\ca m}(\cdot,s,1;\overline{f}),E_\cb(\cdot, 1/2+ir)\rangle E_\cb(z_0, 1/2+ir)\, dr\\
 =\sum_l P_l(m,s) \int_T^{T+1}  \frac{\langle Q_{\ca m}(\cdot,s+W+c_l,1-d_l;\overline{f}),E_\cb(\cdot, 1/2+ir)\rangle}{R_l(1/2+ir,s)}  E_\cb(z_0, 1/2+ir)\, dr \\
 =\sum_l P_l(m,s) \int_T^{T+1} \int_{\F} \frac{Q_{\ca m}(z,s+W+c_l,1-d_l;\overline{f})}{R_l(1/2+ir,s)}  \overline{E_\cb(z, 1/2+ir)} E_\cb(z_0, 1/2+ir)\, d\mu z \, dr.
\endmultline\tag 11.13
$$
The integrand satisfies
$$
\frac{Q_{\ca m}(z,s+W+c_l,1-d_l;\overline{f})}{R_l(1/2+ir,s)}  \overline{E_\cb(z, 1/2+ir)} E_\cb(z_0, 1/2+ir) \ll |r|^{-2W}y_\G(z)^{1/4-W/2}y_\G(z)^{1/2}y_\G(z_0)^{1/2}
$$
by (11.3), (11.4), (11.6) and (11.12). Thus the double integral in (11.13) is absolutely and uniformly convergent and we may interchange the limits of integration to obtain
$$
\sum_l P_l(m,s)  \int_{\F} Q_{\ca m}(z,s+W+c_l,1-d_l;\overline{f})\int_T^{T+1}\frac{\overline{E_\cb(z, 1/2+ir)}}{R_l(1/2+ir,s)}   E_\cb(z_0, 1/2+ir)\, dr \, d\mu z. \tag 11.14
$$
Also
$$
\multline
\int_T^{T+1}\frac{\overline{E_\cb(z, 1/2+ir)}}{R_l(1/2+ir,s)}   E_\cb(z_0, 1/2+ir)\, dr \\
\ll
T^{-2W} \sqrt{ \int_T^{T+1} |E_\cb( z, 1/2+ir)|^2\, dr \cdot \int_T^{T+1} |E_\cb( z_0, 1/2+ir)|^2\, dr}. 
\endmultline
$$
So, with Lemma 8.3 (for $n=0$), (11.14) is bounded by a constant times
$$
\multline
\sum_l |P_l(m,s)|  \int_{\F} |Q_{\ca m}(z,s+W+c_l,1-d_l;\overline{f})| T^{-2W}y_\G(z)^{1/2}T^{6}y_\G(z_0)^{1/2}T^{6} \, d\mu z \\
\ll 
\sum_l |P_l(m,s)|  T^{12-2W} \int_{\F}  y_\G(z)^{3/4-W/2} \, d\mu z \,\cdot  y_\G(z_0)^{1/2}.
\endmultline
$$
This means that, for $W$ chosen large enough,
$$
\int_{-\infty}^{\infty}\langle
Q_{\ca m}(\cdot,s,1;\overline{f}) ,E_\cb(\cdot,1/2+ir)\rangle E_\cb(z,1/2+ir)\,dr
\ll 
y_\G(z)^{1/2}.
$$
Combine this with (11.9) and (11.10) to see that $Q_{\ca m}(z,s,1;\overline{f})$ is analytic for $\Re(s)> 1- \delta_\G$ and bounded by $y_\G(z)^{1/2}$ except for a simple pole at $s=1$ with the stated residue.

\vskip .10in
We leave it to the reader to check the result for $R_0Q_{\ca m}(z,s,1;\overline{f})$ by applying $R_0$ to both sides of (11.1) and using the estimates from Section 8 and Proposition E. Note that $R_0$ eliminates the pole at $s=1$ coming from $\eta_0$, the constant eigenfunction. This completes the proof of Proposition D. \bs

\flushpar {\bf Proof of Proposition E: } We begin with the formula
$$
\overline{f^{(n)}(\g z)}=(-2i)^{-n} \Im(\g z)^{-n-1} \sum_{r=0}^n (-1)^{n-r} \ee(\g, z)^{-2r-2}\binom{n}{r} \frac{(n+1)!}{(r+1)!} L^r\left(y \overline{f(z)}\right)
$$
for $f$ in $S_2(\G)$ and $\g$ in $\G$. This formula may be proved by induction, see \cite{CO}. 
Then, by definition, $I_{\ca(-n)}(z)=f_\ca^{(n)}(z)$ for $-n\leqslant 0$ and
$$
Q_{\ca m}(z,s,-n;\overline{f})=\sum_{\g \in \G_\ca \backslash \G}  \overline{f_\ca^{(n)}(\sa^{-1} \g z)} \Im(\sa^{-1} \g z)^s e(m \sa^{-1} \g z).
$$
So, if we name $\G'=\sa^{-1}\G \sa$ and note that $\sa^{-1}\G_\ca \sa=\G_\ci$, we find
$$
\multline
Q_{\ca m}(\sa z,s,-n;\overline{f})=\sum_{\g' \in \G_\ci \backslash \G'}  
\overline{f_\ca^{(n)}( \g' z)} \Im( \g' z)^s e(m  \g' z)\\
=(-2i)^{-n} \sum_{r=0}^n (-1)^{n-r} \binom{n}{r} \frac{(n+1)!}{(r+1)!}  
L^r\left(y \overline{f_\ca(z)}\right)\sum_{\g' \in \G_\ci \backslash 
\G'}\Im(\g' z)^{s-n-1}  e(m  \g' z) \ee(\g', z)^{-2r-2}\\
=(-2i)^{-n} \sum_{r=0}^n (-1)^{n-r} \binom{n}{r} \frac{(n+1)!}{(r+1)!} 
L^r\left(y \overline{f_\ca(z)}\right) \ee(\sa, z)^{-2r-2} U_{\ca m}(\sa 
z, s-n-1, 2r+2).
\endmultline
$$
Thus, recalling (8.20), (8.21)
$$
\align
 L^r \left(y \overline{f_\ca(z)}\right)  &= L^r \theta_{\sa, -2}\left(y \overline{f(z)}\right)\\
&= \theta_{\sa, -2r-2}L^r \left(y \overline{f(z)}\right)\\
&= L^r \left .\left(y \overline{f(z)}\right)\right|_{\sa z} \ee(\sa, z)^{2r+2}.
\endalign
$$
So we get
$$
Q_{\ca m}(z,s,-n;\overline{f})=(-2i)^{-n} \sum_{r=0}^n (-1)^{n-r} \binom{n}{r} \frac{(n+1)!}{(r+1)!} L^r\left(y \overline{f(z)}\right)  U_{\ca m}( z, s-n-1, 2r+2). \tag 11.15
$$
This identity (11.15) provides the analytic continuation of $Q_{\ca m}(z,s,-n;\overline{f})$ to $\Re(s)>2+n-\delta_\G$ by Propositions B and C.

\vskip .10in
The piece $L^r\left(y \overline{f(z)}\right)$ has exponential decay at every cusp $\cb$ because
$$
\align
 \theta_{\sb, -2r-2} L^r \left(y \overline{f(z)}\right)  &= L^r \left(\theta_{\sb, -2} y \overline{f(z)}\right)\\
&= L^r \left(y \overline{ j(\sb z)^{-2} f(\sb z)}\right)\\
&= L^r \left(y \sum_{n=1}^\infty \overline{ a_\cb (n)e(nz)}\right).
\endalign
$$
Hence 
$$
L^r \left(y \overline{f(z)}\right) \ll y_\G(z)^{r+1} e^{-2\pi y_\G(z)} \tag 11.16
$$
for an implied constant depending on $r, f$ and $\G$. Therefore, with (11.15), (11.16), Propositions B and C, we have
$$
Q_{\ca m}(z,s+n+1,-n;\overline{f})  \ll  e^{-\pi y_\G(z)} 
$$
say, for $\Re(s)>1-\delta_\G$.

\vskip .10in
To show that the same results are true for $R_0 Q_{\ca 
m}(z,s+n+1,-n;\overline{f})$ apply $R_0$ to both sides of (11.15) and 
note that
$$
\multline
R_0\left( L^r\left(y \overline{f(z)}\right)  U_{\ca m}( z, s-n-1, 2r+2)\right) \\
=
\left( R_{-2r-2}L^r\left(y \overline{f(z)}\right)\right)  U_{\ca m}( z, s-n-1, 2r+2)+
 L^r\left(y \overline{f(z)}\right)  R_{2r+2}U_{\ca m}( z, s-n-1, 2r+2)
\\
=
 R_{-2r-2} L^r\left(y \overline{f(z)}\right)  U_{\ca m}( z, s-n-1, 2r+2)\\
+
 L^r\left(y \overline{f(z)}\right)  \left( (s-n+r)U_{\ca m}( z, s-n-1, 2r+4)-4\pi m U_{\ca m}( z, s-n, 2r+4)\right)
\endmultline
$$
by (8.13). This completes the proof of Proposition E. \bs

$$\text{\bf 12. Further Questions}$$

Many natural questions arise:
\roster
\item We have found the dimensions of the spaces $S_k^2(\G)$ and $M_k^2(\G)$ 
focusing on even weight $k$. What are the dimensions for $k$ odd (as in Theorem 2.25 in \cite{Sh} where the dimensions of $S_k(\G)$ and $M_k(\G)$  are given) or for $\G \backslash \H$ compact? 
\item As we saw in $(3.12)$, forms with characters arise naturally in the work of Kleban and Zagier. All the spaces we have discussed may be generalized to arbitrary weights and multiplier systems. 
\item What are the Fourier coefficients of these second-order forms and 
do they have arithmetic or other significance? 
\item
What is the dimension of $R_k^2(\G)$? It is certainly true by Remark 5.5 that 
$$
\dim S_2^2(\G) < \dim R_2^2(\G) \leqslant \dim M_2^2(\G).
$$
\item Is there a natural inner product on the spaces $S_k^2(\G)$ and $M_k^2(\G)$ that respects (3.10) and (3.11)?
\item
A further interesting extension of this work is to higher order forms. We purposely designed our notation in Section 2 wih this in mind. Define the third-order space $S_k^3(\G)$ with the conditions 
{\bf H}, {\bf A(}$S^{2}_k(\G)${\bf )}, {\bf P} and {\bf C}. 
Recursively, $S_k^n(\G)$ should satisfy {\bf H}, {\bf A(}$S^{n-1}_k(\G)${\bf )}, {\bf P} and {\bf C} so that the automorphy condition involves a form of lower order. Similarly for the higher order versions of the other spaces in Section 2. We expect that the methods used in this paper should generalize to counting dimensions of these higher order spaces.
\endroster

\bf Acknowledgements. \rm The authors thank R. Bruggeman and M. Knopp for 
their very useful comments and suggestions.

\end{document}